\documentclass[11pt]{article}
\usepackage[utf8]{inputenc}
\usepackage[dvipsnames]{xcolor}
\usepackage{tikz}
\usetikzlibrary{shapes,arrows,positioning}
\usetikzlibrary{automata,arrows,positioning,calc}
\usepackage{amsmath}
\usepackage{amssymb}
\usepackage{amsthm}
\usepackage{hyperref}
\usepackage{bbm}
\usepackage{blkarray}
\usepackage{empheq}
\usepackage{enumitem}

\usepackage{algorithm}
\usepackage[noend]{algpseudocode}
\usepackage{caption}
\usepackage{subcaption}
\usepackage{float}
\usepackage{stmaryrd}
\usepackage{fancyhdr}
\usepackage{mathabx}
\makeatletter
\newcommand{\pushright}[1]{\ifmeasuring@#1\else\omit\hfill$\displaystyle#1$\fi\ignorespaces}
\newcommand{\pushleft}[1]{\ifmeasuring@#1\else\omit$\displaystyle#1$\hfill\fi\ignorespaces}
\makeatother

\usepackage{amsopn}

\usepackage{bbm}

\usepackage{booktabs}
\renewcommand{\arraystretch}{1.2}

\hypersetup{colorlinks,
citecolor=orange,
citebordercolor=magenta,
linkcolor=magenta
}

\providecommand{\keywords}[1]
{
  \small	
  \textbf{{Keywords.}} #1
}

\numberwithin{equation}{section}

\usepackage{todonotes}

\newtheorem{theorem}{Theorem}[section]
\newtheorem{remark}[theorem]{Remark}
\newtheorem{definition}[theorem]{Definition}
\newtheorem{lemma}[theorem]{Lemma}

\newtheorem{condition}[theorem]{Condition}

\DeclareMathOperator*{\argmin}{arg\,min}

\title{Mean Field Models to Regulate Carbon Emissions in \\Electricity Production}

\author{Ren\'e Carmona\footnote{Department of Operations Research and Financial Engineering,
  Princeton University, 
  Princeton, NJ 08544 
  (\href{mailtorcarmona@princeton.edu}{rcarmona@princeton.edu},
  \href{mailto:gokced@princeton.edu}{gokced@princeton.edu},
  \href{mailto:lauriere@princeton.edu}{lauriere@princeton.edu}).}
\and G\"ok\c ce Dayan{\i}kl{\i}
    \footnotemark[1]
\and Mathieu Lauri\`ere
    \footnotemark[1]
}
\date{}
\allowdisplaybreaks
\begin{document}

\maketitle
\begin{abstract}
The most serious threat to ecosystems is the global climate change fueled by the uncontrolled increase in carbon emissions. In this project, we use mean field control and mean field game models to analyze and inform the decisions of electricity producers on how much renewable sources of production ought to be used in the presence of a carbon tax. The trade-off between higher revenues from production and the negative externality of carbon emissions is quantified for each producer who needs to balance in real time reliance on 
reliable but polluting (fossil fuel) thermal power stations versus investing in and depending upon clean production from
uncertain wind and solar technologies.
We compare the impacts of these decisions in two different scenarios: 1) the producers are competitive and hopefully reach a \textit{Nash Equilibrium}; 2) they cooperate and reach a \textit{Social Optimum}. We first prove that both problems have a unique solution using forward-backward systems of stochastic differential equations. We then illustrate with numerical experiments the producers' behavior in each scenario. We further introduce and analyze the impact of a regulator in control of the carbon tax policy, and we study the resulting Stackelberg equilibrium with the field of producers.
\end{abstract}

\vskip3mm

\vskip 12pt\noindent
\emph{\textbf{Acknowledgments.}}
{The authors were partially supported by NSF DMS-1716673, ARO  W911NF-17-1-0578 and AFOSR FA9550-19-1-0291.}

\section{Introduction} 
\label{sec:intro}
Nowadays, it is widely accepted that the most serious threat to ecosystems is the global warming fueled by the uncontrolled increase in carbon emissions, and for the last twenty-some years, starting with the Kyoto Protocol in 1997, international treaties have sprung out in hope to address this negative externality. The most recent of these treaties is the \textit{Paris Agreement} with $196$ signatories aiming at keeping the increase in temperature below $2 ^\circ$C. Throughout the world, local and federal governments try to disincentivize reliance on polluting means of production by introducing carbon taxes or cap-and-trade programs. In the latter case, regulators put a limit on the allowable quantity of Green House Gas (GHG) emissions, any quantity above this limit having to be covered by emission certificates (allowances) or the payment of a penalty. In the former case, whether they are levied \textit{upstream} or \textit{downstream}, carbon taxes aim at penalizing the use of fossil fuels for their carbon content. The interested reader is referred to \cite{carmona_SICON,carmona_SIREV,carmona_ASCONA} for a review of the state of affairs in the early days of the European Union Emission Trading System, and mathematical treatments of thorough partial equilibrium models for the comparison of realistic implementations of these policies in the electricity sector.

According to the Environmental Protection Agency, electricity production claims the lion share ($25\%$) of the total Greenhouse Gas emissions in the US. \footnote{\url{https://www.epa.gov/ghgemissions/global-greenhouse-gas-emissions-data}} So here, we concentrate on the electricity sector and we propose a model for the analysis of the impact of investments in clean means of production (e.g. solar and wind).  
While a model of the electricity sector should comprise at least three types of agents:  electricity producers,  resellers / retailers, and the end-users, we shall concentrate our modeling effort on the producers. Until the challenging technological problem of electricity storage is resolved at a larger scale, the demand for this commodity remains inelastic, and we shall penalize the producers for not matching the demand, forcing the Independent System Operator (ISO) to rely on costly reserves.
In the following, we shall use the term \textit{renewable} to mean electricity produced from wind turbines or solar panels. Alternatively, we shall use the term \textit{non-renewable} to mean electricity produced by burning fossil fuels like coal, crude oil or natural gas.
We chose this convention for convenience, even if this literary license is not completely accurate.

Individual producers control over time their usage of fossil fuels, and hence, the amount of $CO_2$ emissions they are responsible for. They also control their possible investment in solar or wind production, should they decide to go that route. Notice that while the decision to use fossil fuels changes over time, the investment in solar panels or wind turbines is a one-time decision made at the beginning of the time period under consideration. In our model, producing electricity from  renewable sources involves an initial investment and no extra cost over time since the marginal cost of running these production assets is practically zero (except from maintenance costs and possible subsidies which we ignore here). While the zero cost of production is an attractive feature, it comes with the very high risks due to the difficulties to predict the weather and the uncertainty associated with the high volatility of these predictions. On the other hand, production from traditional power plants is more predictable, the costs depending upon the prices of the fuels and the price put on the $CO_2$ emissions by the regulator. Each producer has to find the right balance between the pros and the cons of the two major means of production we single out in our stylized model. The overarching goal is to decarbonize so as to meet emission targets, harnessing demand-side policies through the establishment of a tax, as well as supply-side resources including wind and solar production technologies.

Our economic model is based on the premises that the individual producers and the regulator have only access to aggregate quantities. Basically, they only have access to the statistical distributions of the productions, emissions, investments, etc of the individual producers.
As a result, we propose two separate frameworks for the individual producers to optimize the mix of renewable and nonrenewable production they should include in their portfolios. We compute and compare the optimal centralized strategies by solving mean field control problems, and the optimal decentralized strategies by solving mean field game problems. Our theoretical analysis relies on the probabilistic approach to construct forward-backward stochastic differential equation (FBSDE) systems for which we show, in both settings, existence and uniqueness of the solutions. Further, we propose a numerical approach to monitor the effect of a carbon tax on the optimal and equilibrium decisions in both cases. Quantifying the differences between the two approaches is reminiscent of what is known as the Price of Anarchy (PoA). 

Among the conclusions drawn from the analysis of our model, we confirm that a carbon tax is an effective incentive for the use of renewables. Also intuitive is the fact that in the absence of a carbon tax, the overall pollution is greater when producers compete than when they cooperate. Less obvious is the fact that cooperating producers will pollute less than when they compete, even if the carbon tax is significant. We also show that stricter regulations tend to reduce the differences between competitive and cooperative equilibria. Further, we argue that the best way for the regulator to encourage producers to match the demand is to incentivize competition over cooperation among the producers.

Mean Field Game (MFG) models appeared simultaneously and independently in the original works of \cite{caines_huang_malhame_2006} and \cite{lasry_lions_2007}. The thrust of these works was to propose a paradigm to overcome the challenges of the search for Nash equilibria in large games by considering models for which the interactions between the players were of a mean field type, and deriving effective equations in the limit when the number of players goes to infinity. Models in which a single player plays a different role from the field of remaining players were introduced and studied under the name of MFGs with major and minor players. In their Stackelberg version, they had a significant impact on problems in economic contract theory. See for example \cite{bensoussan_yam_2016}, \cite{salhab2016},  \cite{posamai_2019}, or \cite{wang_2020}, \cite{elie_2020} or \cite{aurell2020optimal}. Notice that in these models, the major player uses a time dependent control, while in this paper, we shall assume that the regulator uses time independent controls.

Using mean field models for energy applications is very natural. Competition in the oil industry and the impact of the renewable energy competition was analyzed in \cite{gueant_2010} and \cite{sircar_2017}. The early work  \cite{gueant_2010} was extended with the addition of a regulator in \cite{achdou_2020}. In \cite{aid_2020}, optimal entry and exit times for two types of agents, electricity producers using either renewable or nonrenewable energy resources, are analyzed using MFGs. Competition among electricity producers is analyzed in \cite{djehiche_2018} by using Mean Field Type Game where the mean field interactions come through conditional expectation of the electricity price and in \cite{alasseur_2020} by using a model where the interactions enter the electricity spot price. In \cite{huyen_2020} and \cite{elie_2020}, electricity consumers constitute the mean field population and a single electricity producer plays the role of the principal, in contrast to our model where we take the electricity producers as the mean field population and the regulator as the principal.

Mean field models have also been used to model environmental impacts. In \cite{malhame_2017}, a MFG model is proposed to model climate change negotiations among countries interacting through a $CO_2$ emission permit market. Emission certificate markets are also studied in \cite{shrivats_2020} and \cite{zhang_2016}, again without the presence of a regulator.

The paper is structured as follows. In Section~\ref{sec:minormodel}, we introduce the minor players' model and the various equilibrium notions used in the sequel. In Section~\ref{sec:minor_main_theor_res} (resp.  \ref{sec:minor_main_numeric_res}), the main theoretical (resp. numerical) results for the minor players' model are given. In Section~\ref{sec:regulator_model}, we introduce the regulator and define the relevant notions of equilibrium. Finally, we provide numerical results for the combined model with minor players and the regulator in Section~\ref{sec:reg_main_numeric_res} and we summarize our findings in a short Section~\ref{sec:conclusion}.

\section{Mean Field Model for Electricity Producers}\label{sec:minormodel}

\subsection{N-Player Model}

Although we will focus on mean field limits involving an infinite number of players, we start with the description of what the \textbf{$\mathcal{N}$-player version} of  the game would be.  For symmetry reasons, we assume that the total electricity demand is split equally between all the agents, and each agent faces the same demand, say $D_t$ at time $t$. The state of producer $i$ is five-dimensional: instantaneous electricity production $Q^i_t \in \mathbb{R}_+$, instantaneous irradiance $S^i_t \in \mathbb{R}_+$, instantaneous emission level $E^i_t \in \mathbb{R}_+$, cumulative pollution $P^i_t \in \mathbb{R}_+$, and instantaneous nonrenewable energy production $\tilde{N}^i_t \in \mathbb{R}_+$. 
Producer $i$ controls their state by choosing at time $t=0$, their initial investment $R^i_e \in \mathbb{R}$ in renewable production assets (e.g. the number of solar panels they purchase), and at each subsequent time $t$, by choosing the rate of change $N^i_t \in \mathbb{R}$ in nonrenewable energy production.
Notice that $N^i_t$ is time dependent while $R^i_e$ is time independent. This will be a challenging feature of the mathematical analysis of our model.

\begin{remark}
For the sake of definiteness, we use the terminology of solar power production. However, other types of renewable energy can be modelled in a similar way. For example, for wind power, $S^i_t$ would stand for the instantaneous output of a wind farm and $R_e$ would be the corresponding units of initial investment.
\end{remark}

So with these proviso out of the way, we define the time evolution of the state of producer $i$ as:
\begin{equation*}
\label{eq:nplayer_dynamics}
    \begin{alignedat}{2}
    dQ^i_t &= \underbracket{\kappa_{1} N^i_tdt}_{\text{Term~1}} + \underbracket{ \kappa_{2}R_e^i\left( \alpha \cos(\alpha t) dt + \textcolor{Bittersweet}{dS^i_t}\right)}_{\text{Term~2}},  &&\\ 
    \textcolor{Bittersweet}{dS^i_t} &= (\theta - S^i_t) dt + \sigma_0 d\widecheck W^i_t, &&dE^i_t = \delta  N^i_t dt + \sigma_1 d W^i_t,\\ 
    dP^i_t &= E^i_t dt, &&d \tilde{N}^i_t = N^i_t dt.
    \end{alignedat}
\end{equation*}    
The instantaneous electricity production changes depend on the instantaneous nonrenewable energy  usage (given by term 1) and the instantaneous yield from the renewable energy investment (given by term 2). 
This second term includes a seasonality component (sinusoidal term) and a random shock for the variability of the sun irradiance. The form of the seasonality component was chosen for the sake of simplicity. It can easily be extended to several harmonics to include nightly and daily, monthly and yearly effects. In any case, we have 
$
Q_t^i= Q^i_0 + \kappa_1 \tilde{N}^i_t + \kappa_2 R^i_e (\sin(\alpha t) + S^i_t)
$
where $\kappa_1, \kappa_2>0$ are constants that give the efficiency of the production from nonrenewable and renewable energy, respectively. The constant $\alpha>0$ gives the period of the seasonality of the renewable energy.  

We model the idiosyncratic noise terms $S^i_t$ in the renewable productions as independent stationary processes. For the sake of definiteness, we assume that they are Ornstein-Uhlenbeck processes with the same mean $\theta>0$ and volatility $\sigma_0>0$, the $\widecheck W^i$ being independent Wiener processes.

The dynamics of the instantaneous emissions $E^i_t$ have two components: the contribution from the production from nonrenewable energy power plants, and idiosyncratic random shocks with constant volatility  $\sigma_1>0$
given by independent Wiener processes $W^i$, also independent of the $\widecheck W^i$'s. The choice of the constant $\delta$ could include the effects of some abatement measures such as carbon capture, sequestration and the use of filters.

Using the notation $\tilde N^i_t$ for the instantaneous nonrenewable given by
$
\tilde N^i_t = \tilde N^i_0 + \int_0^t N^i_s ds
$,
the expected cost to producer $i$ over the whole period is:
\begin{multline}
\label{eq:nplayer_minorcost} 
    C^{\mathcal{N}}(N^i, R^i_e; \bar Q)
    = 
    \mathbb{E}\Big[\int_0^T \Big[\underbracket{c_{1} |N^i_t|^2}_{\text{Term~1}} +
    \underbracket{p_1 \tilde{N}^i_t}_{\text{Term~2}}+
    \underbracket{c_2|Q^i_t-D_t|^2}_{\text{Term~3}} - \\
    \underbracket{c_3\big(\rho_0 - \rho_1(D_t - \bar Q_t)\big)Q^i_t}_{\text{Term~4}}\Big] dt + 
    \underbracket{\tau|P^i_T|^2}_{\text{Term~5}} + 
    \underbracket{p(R_e^i)\Big]}_{\text{Term~6}},
\end{multline}
where $\bar Q = \sum_{j=1}^{\mathcal{N}} Q^j / \mathcal{N}$ and $p: \mathbb{R}_+ \mapsto \mathbb{R}_+$ is the price function for the investment in renewable energy.

Term~1 with $c_1>0$, is a penalty (ie. delay cost) for attempting to ramp up and down nonrenewable energy power plants too quickly. 
Term~2 represents the costs of the fossil fuels used in nonrenewable power plants. The constant $p_1>0$ can be understood as the average cost of one unit of fossil fuel.
In lieu of storage which is not included in our models because of its scarcity, Term~3 with $c_3>0$, imposes a penalty on producers for not matching the demand  and forcing the system operator to use costly reserves.
Term~4 represents the revenues from electricity production, $\big(\rho_0 + \rho_1 (D_t - \bar Q_t)\big)$ being the inverse demand function which is assumed to be linear in excess demand or supply. Here $\rho_0$and  $\rho_1$ are strictly positive constants. It captures the fact that the price increases if there is excess demand, and it decreases if there is excess supply. We assume that the producers are selling what they produce. This term introduces the mean field interactions into the model. 
Term~5 gives the carbon tax levied by the regulator. We emphasize its role by assuming it is proportional to the square of the terminal pollution. 
Term~6 is the total cost related to the initial investment in renewable electricity production including the price of the solar panels and the cost of the land used.

\subsection{The Mean Field Model}

In discussing the mean field regime of the model, we focus on a \textit{representative} producer interacting with the field of the other producers, so we drop the superscript  $i$ and the dynamics equations become:
\begin{equation}
\label{eq:minordynamics} 
\begin{aligned}
    dQ_t &= \kappa_{1} N_tdt +  \kappa_{2}R_e\big( \alpha \cos(\alpha t) dt + \hskip-3mm &&\textcolor{Bittersweet}{(\theta- S_t) dt + \sigma_0 d\widecheck W_t}\big), \\
    \textcolor{Bittersweet}{dS_t} &= (\theta - S_t) dt + \sigma_0 d\widecheck W_t, && dE_t = \delta  N_t dt + \sigma_1 d W_t,\\
    dP_t &= E_t dt, &&d \tilde{N}_t = N_t dt,
\end{aligned}
\end{equation}
where $W$ and $\widecheck W$ are independent Wiener processes.
Accordingly, the total expected cost becomes:
\begin{multline}
\label{eq:minorcost}  
    C(N, R_e; \bar Q) = 
    \mathbb{E}\Big[\int_0^T \Big[c_{1} |N_t|^2 +
    p_1 \tilde{N}_t+
    c_2|Q_t-D_t|^2 - \\
    c_3\big(\rho_0 + \rho_1(D_t- \bar Q_t)\big)Q_t\Big] dt + 
    \tau|P_T|^2 + 
    p(R_e)\Big],
\end{multline}
where $\bar Q_t=\mathbb{E}[Q_t]$. We shall sometimes use the notation $\bar Q_t(N,R_e)$ to emphasize the fact that the expectation is computed under the state dynamics controlled by the admissible control $(N,R_e)$.

\subsection{Equilibrium Notions}

We consider two different models: mean field game (MFG) and mean field control (MFC). In the mean field game model, producers behave competitively and minimize their total expected costs (search for their best responses) given the other players' decisions. A Nash equilibrium is then characterized as a fixed point of the best response map so defined.  In the sequel, we restrict our attention to admissible strategies $(N,R_e)$ such that $\mathbb{E}[\int_0^T|N_t|^2 dt] <+\infty$ and $R_e \in \mathbb{R}_+$. 

\begin{definition}[\textit{MFG Nash Equilibrium}] An admissible strategy and mean field flow tuple, $(\hat N, \hat R_e, \bar Q)$,  is called an \textit{MFG Nash equilibrium} for any admissible $(N, R_e)$, we have:
$$
C\Big(\textcolor{Bittersweet}{(N, R_e)}; \bar Q\Big) \geq C\Big((\hat N, \hat R_e); \bar Q\Big),
$$
and $\bar Q = \bar Q(\hat N, \hat R_e)$.
\end{definition}
In the mean field control case, we assume that the producers cooperate and leave the choice of the control to a social planner minimizing the total expected cost as defined in \eqref{eq:minorcost}. In the realistic setup, the producers can be thought as the production facilities of a monopolistic electricity production firm and the social planner's decisions refer to the decisions taken by the headquarter. In this case, if one player changes their behavior, every player changes in the same way, and the mean field is affected. The problem is now an optimal control problem. 

\begin{definition}[\textit{Social Planner's MFC Optimum}] An admissible strategy and mean field flow tuple, $(\hat N,\hat R_e, \bar Q)$,  is called an \textit{MFC optimum} if for any admissible $(N, R_e)$, we have:
$$    
C\Big(\textcolor{Bittersweet}{(N, R_e)}; \bar Q(\textcolor{Bittersweet}{N, R_e})\Big) \geq C\Big((\hat N, \hat R_e); \bar Q\Big),
$$
and $\bar Q = \bar Q(\hat N, \hat R_e)$.
\end{definition}

\section{Main Theoretical Results}
\label{sec:minor_main_theor_res}

In this section, the following forward backward stochastic differential equation system (FBSDE) is going to be of interest:
\begin{align}
\label{eq:fbsde}
        dQ_t &= -\frac{\kappa_{1}}{2c_1} (Y_t^1\kappa_1+Y_t^3\delta+Y_t^5)dt \nonumber\\
        &\qquad+  \kappa_{2} (p^{\prime})^{-1} \Big(-\mathbb{E}\Big[\int_0^T \kappa_2Y_t^1\Big(\alpha \cos(\alpha t) + (\theta {-S_t}) \Big) dt \Big]\Big)\nonumber\\
        &\pushright{\times \left( \alpha \cos(\alpha t) dt + (\theta- S_t) dt + \sigma_0 d\widecheck W_t\right),}  &Q_0 &= q_0 \nonumber\\
        dS_t &= (\theta - S_t) dt + \sigma_0 d\widecheck W_t,  &S_0 &= \theta \nonumber\\ 
        dE_t &= -\frac{\delta}{2c_1}  (Y_t^1\kappa_1+Y_t^3\delta+Y_t^5) dt + \sigma_1 d W_t,
         &E_0 &= e_0\nonumber\\
        dP_t &= E_t dt,  &P_0 &= p_0
        \nonumber\\
        d \tilde{N}_t &= \frac{1}{2c_1} (Y_t^1\kappa_1+Y_t^3\delta+Y_t^5) dt, \qquad &\tilde{N}_0 &= \tilde{n}_0\nonumber\\
        dY^1_t &= \Big(-{2c_2({Q}_t-D_t)} +
        {c_3\big(\rho_0 + \rho_1(D_t- \bar Q_t)\big)}\Big)dt + & &\nonumber\\
        & \qquad\qquad\qquad\qquad\qquad\qquad\qquad\qquad Z_t^{1,1}d\widecheck{W}_t + Z_t^{1,2}dW_t    &Y_T^1 &= 0 \nonumber\\
        dY^2_t &= \Big(\kappa_2 (p^{\prime})^{-1} \Big(-\mathbb{E}\Big[\int_0^T \kappa_2Y_t^1\Big(\alpha \cos(\alpha t) + (\theta {-S_t}) \Big) dt \Big]\Big) Y_t^1 + 
        \nonumber\\
        &\pushright{Y_t^2\Big)dt 
        + Z_t^{2,1}d\widecheck{W}_t + Z_t^{2,2}dW_t, }    &Y_T^2 &= 0\nonumber\\
        dY^3_t &= -Y_t^4dt + Z_t^{3,1}d\widecheck{W}_t + Z_t^{3,2}dW_t,    &Y_T^3 &= 0   \nonumber\\
        dY^4_t &= Z_t^{4,1}d\widecheck{W}_t + Z_t^{4,2}dW_t, &Y_T^4 &= 2\tau\widehat{P}_T \nonumber\\
        dY^5_t &= -p_1dt + Z_t^{5,1}d\widecheck{W}_t + Z_t^{5,2}dW_t,   &Y_T^5 &= 0.
\end{align}

\begin{theorem}
\label{theorem:fbsde_mfg}
    $(\hat N_t,\hat R_e, \bar Q)$ is a Nash equilibrium if and only if $(\hat N, \hat R_e)$ is given by:
\begin{equation}
    \begin{aligned}
    \label{eq:fbsde_opt_cond}
        \hat N_t &= -\dfrac{Y_t^1\kappa_1+Y_t^3\delta+Y_t^5}{2c_1}, \quad t \in [0,T]\quad\text{and}\quad\\
        \hat R_e&=(p^{\prime})^{-1} \Big(-\mathbb{E}\Big[\int_0^T \kappa_2Y_t^1\Big(\alpha \cos(\alpha t) + (\theta {-S_t}) \Big) dt \Big]\Big),
    \end{aligned}
\end{equation} 
    where $(Q,S,E,P,\tilde{N},Y^1,Y^2,Y^3,Y^4,Y^5)$ is a solution to the FBSDE given in \eqref{eq:fbsde}.\footnote{Here,
$(p^\prime)^{-1}(\cdot)$ refers to the inverse of the first derivative of the function $p(\cdot)$.}
\end{theorem}

\begin{condition}
\label{cond:cond_on_p}
    \begin{enumerate}[label={(\roman*)}]
        \item $p$ is convex.
        \item $(p^{\prime})^{-1}$ is bounded  i.e. $(p^{\prime})^{-1}:\mathbb{R}\mapsto[0, R_e^{\text{max}}]$, continuous and monotone.
    \end{enumerate}
\end{condition}

\begin{theorem}
\label{theorem:fbsde_existence_mfg}
     Assume Condition~\ref{cond:cond_on_p} holds, then there exists a unique Nash Equilibrium mean field flow $\bar Q$.
\end{theorem}

\begin{theorem}
\label{theorem:fbsde_mfc}
$(\hat N, \hat R_e)$ is an MFC optimum if and only if $(\hat N, \hat R_e)$ is given by \eqref{eq:fbsde_opt_cond} where  $(Q,S,E,P,\tilde{N},$ $Y^1,Y^2,Y^3,Y^4,Y^5)$ is a solution to the FBSDE given in \eqref{eq:fbsde} where the equation for $(Y^1_t)_t$ is replaced by 
\begin{equation}
\label{eq:mfc_y1}
dY^1_t = \Big(-{2c_2({Q}_t-D_t)} +
        {c_3\big(\rho_0 + \rho_1(D_t- \textcolor{Bittersweet}{2\bar Q_t})\big)}\Big)dt + Z_t^{1,1}d\widecheck{W}_t + Z_t^{1,2}dW_t.    
\end{equation}
\end{theorem}

\begin{theorem}
\label{theorem:fbsde_existence_mfc}
Assume Condition~\ref{cond:cond_on_p} holds, then there exists a unique mean field control optimum flow $\bar Q$.
\end{theorem}

\section{Numerical Approach}
\label{sec:minor_main_numeric_res}

For numerical purposes, given the technical challenges posed by the solution of the large FBSDE in~\ref{eq:fbsde} with the existence of time dependent and independent controls,
we implement an analytic approach for which we give the details below. For this reason, we first notice that:
\begin{equation*}
    \inf_{(N_t)_t, R_e} C(N, R_e; \bar Q) = \inf_{R_e} \inf_{(N_t)_t} C(N, R_e; \bar Q),
\end{equation*}
and we assume that $R_e$ is fixed in a first analysis.
Next, we rewrite the model in matrix form using $X_t:=[Q_t\quad S_t\quad E_t\quad P_t\quad \tilde{N}_t]^{\top}$ as $5$-dimensional state process at time $t$, and rewrite the optimization problem as: 
\begin{equation}
\begin{aligned}
\label{eq:minorcost_matrix}
 \inf_{(N_t)_t} \tilde{C}\Big(N; R_e, \bar X\Big) =&
\inf_{(N_t)_t} \mathbb{E}\Bigg[\int_0^{T} \Big[\frac{R}{2} |N_t|^2 + H^{\top}_t X_t + \bar{X}_t^{\top}F X_t  + X_t^{\top} G X_t + J_t\Big] dt\\&\pushright{+ X^{\top}_T S_T X_T + p(R_e)}
 \Bigg]
 \end{aligned}
\end{equation}
\begin{equation*}
    dX_t = \Big(A X_t + B \cdot N_t + C_t \Big) dt + \Sigma d\widetilde{W}_t
\end{equation*}
where $R$ and $J_t$ are the scalars given by $R = 2c_1$ and $J_t = c_2 D_t^2$ and: 
{\small\begin{equation*}\arraycolsep3pt 
    H_t =  \begin{bmatrix}
        -(2c_2+c_3\rho_1)D_t-c_3 \rho_0\\
        0\\
        0\\
        0\\
        p_1
    \end{bmatrix},
    F = \begin{bmatrix}
        c_3 \rho_1& 0 & 0 & 0 & 0\\
        0 & 0 & 0 & 0 & 0\\
        0 & 0 & 0 & 0 & 0\\
        0 & 0 & 0 & 0 & 0\\
        0 & 0 & 0 & 0 & 0
        \end{bmatrix},
    G = \begin{bmatrix}
        c_2 & 0 & 0 & 0 & 0\\
        0 & 0 & 0 & 0 & 0\\
        0 & 0 & 0 & 0 & 0\\
        0 & 0 & 0 & 0 & 0\\
        0 & 0 & 0 & 0 & 0
        \end{bmatrix},
    S_T = \begin{bmatrix}
        0 & 0 & 0 & 0 & 0\\
        0 & 0 & 0 & 0 & 0\\
        0 & 0 & 0 & 0 & 0\\
        0 & 0 & 0 & \tau & 0\\
        0 & 0 & 0 & 0 & 0
        \end{bmatrix},
\end{equation*}
\vskip-2mm
\begin{equation*}\arraycolsep3pt 
    A = \begin{bmatrix}
        0 & -\kappa_2 R_e & 0 & 0 & 0\\
        0 & -1 & 0 & 0 & 0\\
        0 & 0 & 0 & 0 & 0\\
        0 & 0 & 1 & 0 & 0\\
        0 & 0 & 0 & 0 & 0
        \end{bmatrix},
    B =  \begin{bmatrix}
        \kappa_1\\
        0\\
        \delta\\
        0\\
        1
    \end{bmatrix}  ,
    C_t =  \begin{bmatrix}
        \kappa_2 R_e\Big( \alpha cos(\alpha t) + \theta \Big)\\
        \theta\\
        0\\
        0\\
        0
    \end{bmatrix},
    \Sigma = \begin{bmatrix}
        \kappa_2 R_e \sigma_0 & 0\\
        \sigma_0 & 0\\
        0 & \sigma_1\\
        0 & 0\\
        0 & 0
    \end{bmatrix}.
\end{equation*}}
Furthermore, we define $\widetilde{W}_t$ and $a$ as:
{\small\begin{equation*}\arraycolsep3pt 
    \widetilde{W}_t = \begin{bmatrix}
        \widecheck W_t\\
        W_t
        \end{bmatrix},\qquad
    a = \frac{1}{2}\Sigma \Sigma^{\top} = \frac{1}{2} \begin{bmatrix}
        (\kappa_2 R_e \sigma_0)^2 & \kappa_2 R_e \sigma_0^2 & 0 & 0 & 0\\
        \kappa_2 R_e \sigma_0^2 &  \sigma_0^2 & 0 & 0 & 0\\ 
        0 & 0 & \sigma_1^2 & 0& 0\\ 
        0 & 0 & 0 & 0 & 0\\
        0 & 0 & 0 & 0 & 0
        \end{bmatrix},
\end{equation*}}
and the value function $u(t,X)$ as:
\begin{equation}
\begin{aligned}
u(t, X)=
\inf_{(N_s)_s} \mathbb{E}\Bigg[\int_t^{T} \Big[\frac{R}{2} |N_s|^2 + H^{\top}_s X_s + \bar{X}_s^{\top}F X_s  + X_s^{\top} G X_s  + J_s\Big] ds\\ \pushright{ + X^{\top}_T S_T X_T + p(R_e)\Big| X_t= X \Bigg]}.
\end{aligned}
\end{equation}

\begin{lemma}[ODE System for the MFG] 
\label{lem:ode_mfg}
For $R_e$ fixed, if there exists a function $t\mapsto (\eta_t, r_t, \bar X_t)$ solving the following system of Ordinary Differential Equations (ODEs):
\begin{subequations}
\label{eq:mfg_ode}
  \begin{empheq}[left=\empheqlbrace]{align}
     &\dfrac{d{\eta_t}}{dt} -{\eta_t} BR^{-1}B^{\top} {\eta_t} + A^{\top} {\eta_t} +{\eta_t} A +2G = 0,  &&{\eta_T} = 2S_T\label{eq:mfg_mfc_eta}\\
    &-\dfrac{d{r_t}}{dt} = \left(A^{\top} -{\eta_t} B R^{-1} B^{\top}\right) {r_t} + {\eta_t} C_t + H_t + F^{\top} \bar{X}_t,  &&{r_T} = 0\label{eq:mfg_r}\\
    &\dfrac{d\bar{X}_t}{dt} = (A-B R^{-1} B^{\top}{\eta_t})\bar{X}_t - B R^{-1} B^{\top}{r_t}+C_t,  &&\bar{X}_0 =\bar{x}_0\label{eq:mfg_mfc_xbar}
  \end{empheq}
\end{subequations}
and if $s_0$ is given by: 
\begin{equation}
\label{eq:mfg_mfc_s}
    s_0 = p(R_e) + \int_0^{T} \Big(tr(a{\eta_t}) -\frac{1}{2} {r_t}^{T} B R^{-1} B^{\top} {r_t} +C_t^{\top} {r_t}+ J_t\Big)dt,
\end{equation}
then $\hat N_t(R_e) = -R^{-1}B^{\top}(\eta_t X_t + r_t)$ is the MFG equilibrium given  $R_e$ fixed, and the expected cost to the representative producer in this equilibrium is:
\begin{equation}
    \label{eq:mfg_cost_eq}
    \inf_{N=(N_t)_t} \tilde{C}^{MFG}\Big(N; R_e, \bar X \Big) = \frac{1}{2} \left( Var(\sqrt{\eta_0} X_0) + \mathbb{E}[\sqrt{\eta_0} X_0]^2 \right) +\bar X_0^{\top} r_0 + s_0.
\end{equation}
\end{lemma}

\begin{theorem}\label{the:exist_uniq_mfg}
For $R_e$ fixed, if $T$ is small enough, there exists a unique MFG equilibrium.
\end{theorem}

\begin{lemma}[MFC ODE System]
\label{lem:ode_mfc}
Given $R_e$, if there exists a function $t\mapsto (\eta_t, r_t, \bar X_t)$ solving the ODE system \eqref{eq:mfg_ode} with the equation \eqref{eq:mfg_r} replaced by:
\begin{equation}
\label{eq:mfc_ode}
    -\frac{d{r_t}}{dt} = \left(A^{\top} - {\eta_t} B R^{-1} B^{\top}\right) {r_t} + {\eta_t} C_t + H_t + F^{\top} {\bar{X}_t} + \textcolor{Bittersweet}{F {\bar{X}_t}}, \qquad {r_T} = 0
\end{equation}
and the same $s_0$ given by \eqref{eq:mfg_mfc_s}, then $N^*_t(R_e) = -R^{-1}B^{\top}(\eta_t X_t + r_t)$ is an optimum for the MFC problem given $R_e$, and the minimal expected cost is 
\begin{equation}
\begin{aligned}
    \label{eq:mfc_cost_eq}
    \inf_{(N_t)_t} \tilde{C}^{MFC}\Big(N; R_e \Big) &= \frac{1}{2} \left( Var(\sqrt{\eta_0} X_0) + \mathbb{E}[\sqrt{\eta_0} X_0]^2 \right) +\bar X_0^{\top} r_0 + s_0\\
    &\pushright{\textcolor{Bittersweet}{-\int_0^T \bar X_t^{\top}F \bar X_t dt}}.
\end{aligned}
\end{equation}

\end{lemma}

\begin{theorem}\label{the:exist_uniq_mfc}
For $R_e$ fixed, if $T$ is small enough, there exists a unique MFC optimum. 
\end{theorem}

Numerically, we search for the $R_e$ and the corresponding equilibrium $N=(N_t)_{t}$ that minimizes the cost of the minor players by using the ODE systems given in \eqref{eq:mfg_ode} and \eqref{eq:mfc_ode}. 

As emphasized earlier, the main difference between MFC and MFG is whether the mean field is affected by the decision of the representative producer (MFC), or taken to be fixed (MFG). This difference translates into the addition of a fixed point argument in the MFG case. For pedagogical reasons, we first discuss the MFC case, then the MFG. After solving the Riccati equation which is the same in both cases, we solve the coupled ODE system directly in the MFC case in order to find the mean field; on the other hand, notice that in the MFG case, the ODEs are decoupled since the mean field is assumed to be fixed in each iteration of the fixed point algorithm.

\subsection{Mean Field Control Algorithm}
In order to solve the system of MFC coupled ODEs for $(\bar{X}_t)_t$ and $r_t$ given by equations \eqref{eq:mfg_mfc_xbar} and \eqref{eq:mfc_ode}, we discretize the time with uniform step size $\Delta t$ and solve the following linear equation:
{\small\begin{equation}
\label{eq:mfc_linearsystem}
    \begin{bmatrix}
        \bar{X}\\
        r
    \end{bmatrix}    =
    M \begin{bmatrix}
        \bar{X}\\
        r
    \end{bmatrix} + K,
\end{equation}}
where $\bar{X} = [\bar{X}_0, \bar{X}_{\Delta t} ,\bar{X}_{2\Delta t}, \dots, \bar{X}_T]^{\top}$, $r = [r_0, r_{\Delta t} ,r_{2\Delta t}, \dots, r_T]^{\top}$.

\begin{algorithm}[H]
\caption{\small Computation of the Mean Field Control Cost over $(N_t)_t$ given $R_e$} {\small
\begin{algorithmic}[1]
\Function{\texttt{Optim-MFC-N}}{$R_e$}
\vskip2mm
    \State Calculate $(\eta_t)_t$ by solving the Riccati Equation in \eqref{eq:mfg_mfc_eta}
    \vskip1mm
    \State Solve the coupled $(\bar{X}_t)_t$ and $(r_t)_t$ linear system \eqref{eq:mfc_linearsystem} given $(\eta_t)_t$ and $R_e$
    \vskip1mm
    \State Calculate $s$ given $R_e$, $(r_t)_t$ and $(\eta_t)_t$ using the equation in \eqref{eq:mfg_mfc_s}
    \vskip1mm
    \State Calculate the expected cost associated with $R_e$, $\hat c:=\inf \tilde{C}^{MFC}(N;R_e, \bar X)$ using \eqref{eq:mfc_cost_eq}
    \vskip2mm
    \State \Return ($\hat c, \bar X$)
    \vskip3mm
\EndFunction    
\end{algorithmic}}
\end{algorithm}

\begin{algorithm}[H]
\caption{\small Search for a Social Optimum} 
{\small
\begin{algorithmic}[1]
\Function {\texttt{SocialOpt}}{} 
		\vskip2mm
		\State Search for the optimal $\hat R_e$ where the optimal cost $R_e \rightarrow c(R_e)$ and optimal mean field $R_e \rightarrow \bar X(R_e)$ are computed by \texttt{Optim-MFC-N}
		\vskip1mm
		\State Let $\hat c = c( \hat R_e )$ and $\hat{\bar X}  = \bar X( \hat R_e )$ 
		\vskip2mm
	\State \Return $(\hat c, \hat R_e, \hat{\bar{X}})$
	\vskip3mm
\EndFunction
\end{algorithmic}}
\end{algorithm}

\subsection{Mean Field Game Algorithm}
In Mean Field Game case, since in each iteration it is assumed that the $(\bar{X}_t)_t$ is fixed, the ODE for $(r_t)_t$ in equation \eqref{eq:mfg_r} can be solved directly by using the following linear equation after we discretize  time:
{\small\begin{equation}
\label{eq:mfg_linearsystem_r}
    r    = M_r r+ K_r,
\end{equation}}
where $r = [r_0, r_{\Delta t} ,r_{2\Delta t}, \dots, r_T]^{\top}$.
Then with this $(r_t)_t$, the time discretization of $(\bar{X}_t)_t$ with dynamics given by the equation \eqref{eq:mfg_mfc_xbar} can be written as:
{\small\begin{equation}
\label{eq:mfg_linearsystem_xbar}
        \bar{X}   =
    M_{\bar{X}} 
        \bar{X} + K_{\bar{X}},
\end{equation}}
where $\bar{X} = [\bar{X}_0, \bar{X}_{\Delta t} ,\bar{X}_{2\Delta t}, \dots, \bar{X}_T]^{\top}$.
The numerical algorithms to find the Mean Field Control and Game Equilibria are given in detail in the following sections.

\begin{algorithm}[H]
\caption{\small Computation of the Expected Cost over $(N_t)_t$ given $R_e$, $(\bar X_t)_t$} 
{\small
\begin{algorithmic}[1]
\Function {\texttt{Optim-MFG-N}}{$R_e, (\bar X_t)_t$}
    \vskip2mm		
    \State Calculate $(\eta_t)_t$ by solving Riccati Equation in \eqref{eq:mfg_mfc_eta}
    \vskip1mm
    \State Solve the linear system for $(r_t)_t$ in \eqref{eq:mfg_linearsystem_r} \textcolor{Bittersweet}{given $(\bar{X}_t)_t$}, $R_e$ and $(\eta_t)_t$
    \vskip1mm
    \State Calculate $s$ given $R_e$, $(r_t)_t$ and $(\eta_t)_t$ using the equation in \eqref{eq:mfg_mfc_s}
    \vskip1mm
    \State Calculate the cost associated with $R_e$ \textcolor{Bittersweet}{and $(\bar{X}_t)_t$}, $\hat c:=\inf \tilde{C}^{MFG}(N; R_e, \bar X)$ using \eqref{eq:mfg_cost_eq}    
    \vskip2mm
	\State \Return $\hat c$
	\vskip3mm
\EndFunction
\end{algorithmic}}
\end{algorithm}

\begin{algorithm}[H]
\caption{\small Search for a Nash Equilibrium}
{\small
\begin{algorithmic}[1]
\Function {\texttt{NashEq}}{}
		\vskip2mm
		\State Initialize $(\bar{X}^0_t)_t$
		\vskip2mm
		\While{$||\bar{X}^k - \bar{X}^{k-1}|| > \epsilon$}
		\vskip2mm
        \State Search for the optimal $\hat R_e$ given $\bar{X}^k$ where the optimal cost $(R_e, \bar{X}^{\textcolor{Bittersweet}{k}}) \rightarrow c^{\textcolor{Bittersweet}{k}}(R_e, \bar{X}^{\textcolor{Bittersweet}{k}})$ is computed by \texttt{Optim-MFG-N}
		\vskip1mm
		\State Let $\hat R_e^{\textcolor{Bittersweet}{k}} = \argmin_{R_e} c^{\textcolor{Bittersweet}{k}}(R_e, \bar{X}^{\textcolor{Bittersweet}{k}})$ 
		\vskip1mm
		\State Compute $(\bar{X}^{\textcolor{Bittersweet}{{k+1}}}_t)_t$ given $\hat R_e^k,(\bar{X_t}^{\textcolor{Bittersweet}{k}})_t$ by solving the linear equation \eqref{eq:mfg_linearsystem_xbar} 
		\EndWhile
	    \vskip2mm
	\State Let $\hat R_e = \hat R^{\textcolor{Bittersweet}{k}}_e$, $\hat c = c^{\textcolor{Bittersweet}{k}}( \hat R_e )$ and $\hat{\bar{X}} = \bar{X}^{\textcolor{Bittersweet}{k}}$
    \vskip2mm
	\State \Return $(\hat c, \hat R_e, \hat{\bar{X}})$
	\vskip3mm
\EndFunction
\end{algorithmic}}
\end{algorithm}

\subsection{Numerical Experiments}
\label{subsec:minor_param}
In the numerical experiments reported below, we use the following parameter values:
\vskip5mm
\begin{center}
{\small
{\renewcommand{\arraystretch}{1.2}
\begin{tabular}{ l|c|l } 
 \hline
$ p_1 = 7/\Delta_t$ (dollar/time)& $ \rho_0=40/\Delta_t$ &  $\theta=5 $\\
$ p_2=10^4$, $p_3=10^{-10}$ &$ \rho_1=0.1/\Delta_t$ &  $T =20$ years, $\Delta_t= 10$ days \\
$ c_1=10^{-4}$ (dollar/$10^3$ cu ft$^2$)& $\alpha=40\pi $& $R_e=[0,5\times10^3]$ (10,000 dollars)\\
$ c_3=1$ (dollar/$10^3$ cu ft)& $\delta=0.15$   & $ D_t =2\times 10^4 - 5\times10^2 \cos(80\pi \Delta_t)$\\
$\kappa_1 =0.13$ (MWh/$10^3$ cu ft)&$\sigma_0=0.01 $  &$\bar{X}_0 = [0,\theta,0,0,0]$\\ 
$\kappa_2=0.1$ (MWh/$10^3$ dollars)&$\sigma_1=0.01 $ &$Var[X_0] = [0,0.1,0,0,0]$\\
 \hline
\end{tabular}}}
\end{center}
\vskip5mm
Furthermore, we assume that $p(R_e)=p_2 R_e-p_3\sqrt{R_e(R_e^{\max}-R_e)}+\epsilon$ where $p_2, p_3$ are positive constants and $\epsilon>0$ is a small constant that ensures the nonnegativity of the price of the units of the renewable energy investment. \footnote{We note that this function satisfies the assumptions that are necessary for existence and uniqueness. The fact that $p(\cdot)$ is convex can be justified by the increasing unit costs that comes from the search of land that is large enough to construct the solar panels on. However, in the numerical application, we take the $p_3$ and $\epsilon$ small to have a function that is nearly ``linear".}

We focus on natural gas as the source of nonrenewable energy. In our numerical experiments, we ignore the effect of the COVID-19 pandemic, and we run simulations for $20$ years starting from March 2020. For the cost of solar power, we use the current assumption that a $1$MW solar farm needs a $1$M\$ investment\footnote{\url{https://news.energysage.com/solar-farms-start-one/}}, and use daily peak sun hours data to compute daily average production from solar panels. We assume that on average, peak hours last approximately 5 hours in the US, and we infer that a solar farm built with  an initial investment of $\$ 10,000$ generates on average $0.5$MWh in $10$ days. We choose $\alpha$ to take into account the seasonality, since sun exposure levels are maximum during summers, minimum during winters. Therefore, we infer that one unit investment of $R_e$ corresponds to $\approx$\$10,000 and on average it generates $\kappa_2(\theta \pm 1)=0.1(5\pm1)$MWh electricity in $10$ days.

According to the data provided by the U.S. Energy Information Administration (EIA), in 2018, $1,365,822$ Million KWh were produced by the electric sector by using 10,215 billion cubic feet of natural gas\footnote{\url{https://www.eia.gov/totalenergy/data/monthly/pdf/sec7.pdf}, (Table 7.2b \& 7.3b)}. Therefore, we assume that $1000$ cubic feet of natural gas produce approximately $ 0.13$MWh. Again, according to the data provided by EIA, typical natural gas power plants produce $0.91$ pound of carbon emission per kWh electricity generation.\footnote{\url{https://www.eia.gov/tools/faqs/faq.php?id=74&t=11}} Since $1000$ cubic feet of natural gas produce around $0.13$ MWh, we conclude that around $122$ pounds of carbon is emitted with the use of natural gas. Moreover, if coal or other fossil fuels were used, this carbon pollution would have to be much more than doubled. Taking this into account  and converting to metric tone we end up with $\delta\approx0.15$.  

We assume that the average demand of electricity for each plant is around the capacity of the plants. According to EIA\footnote{\url{https://www.eia.gov/electricity/annual/archive/pdf/epa_2018.pdf}, (Table 4.3)}, the average daily capacity of a natural gas plant is around $2166.5$ MWh in 2018. Furthermore, the monthly seasonal component is found by using the monthly residential electricity consumption in 2018 data given by EIA\footnote{\url{https://www.eia.gov/totalenergy/data/monthly/pdf/sec7.pdf}, (Table 7.6)}. Therefore, $10$ day demand is taken sinusoidal to show the seasonality around $20,000$MWh.
According to the data provided by EIA\footnote{\url{https://www.eia.gov/electricity/annual/html/epa_08_04.html}, \\ \url{https://www.eia.gov/dnav/ng/ng_pri_sum_a_EPG0_PEU_DMcf_a.htm}} nonrenewable energy has 40\% of its fuel cost as the operation and maintenance costs on top of the fuel cost in 2018 and the price of 1000 cu ft natural gas can be assumed \$5. Therefore, we take $p_1=\$7$. Finally by using the daaata given by EIA\footnote{\url{https://www.eia.gov/todayinenergy/detail.php?id=37912}}, we see that the average price of wholesale electricity is around \$40 per MWh, therefore we take $\rho_0=40$.

\subsubsection{Price of Anarchy (PoA) Analysis}

\begin{figure}[H]
\centering
\begin{subfigure}{.44\textwidth}
    \includegraphics[width=1\linewidth]{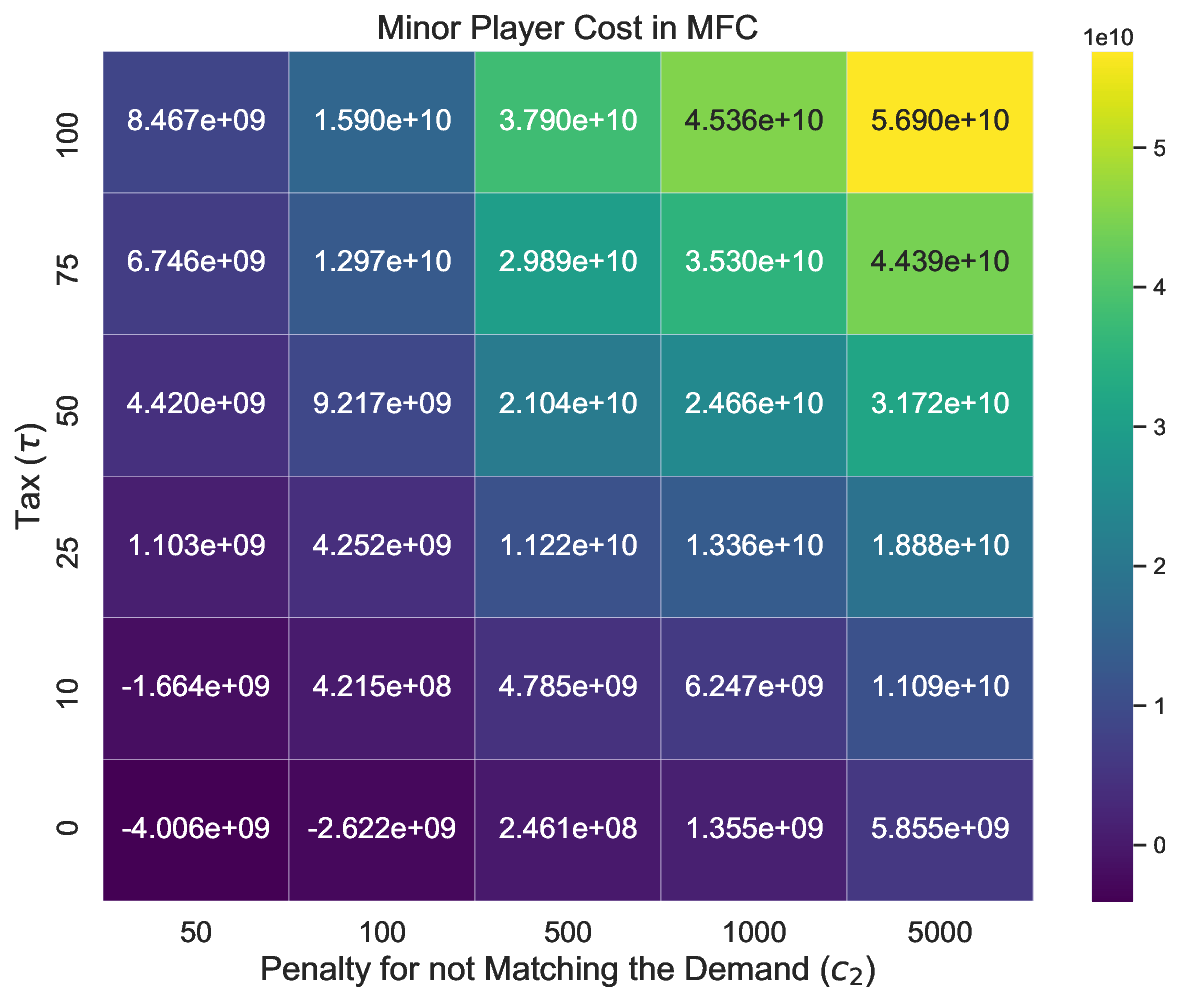}
\end{subfigure}
\begin{subfigure}{.44\textwidth}
    \includegraphics[width=1\linewidth]{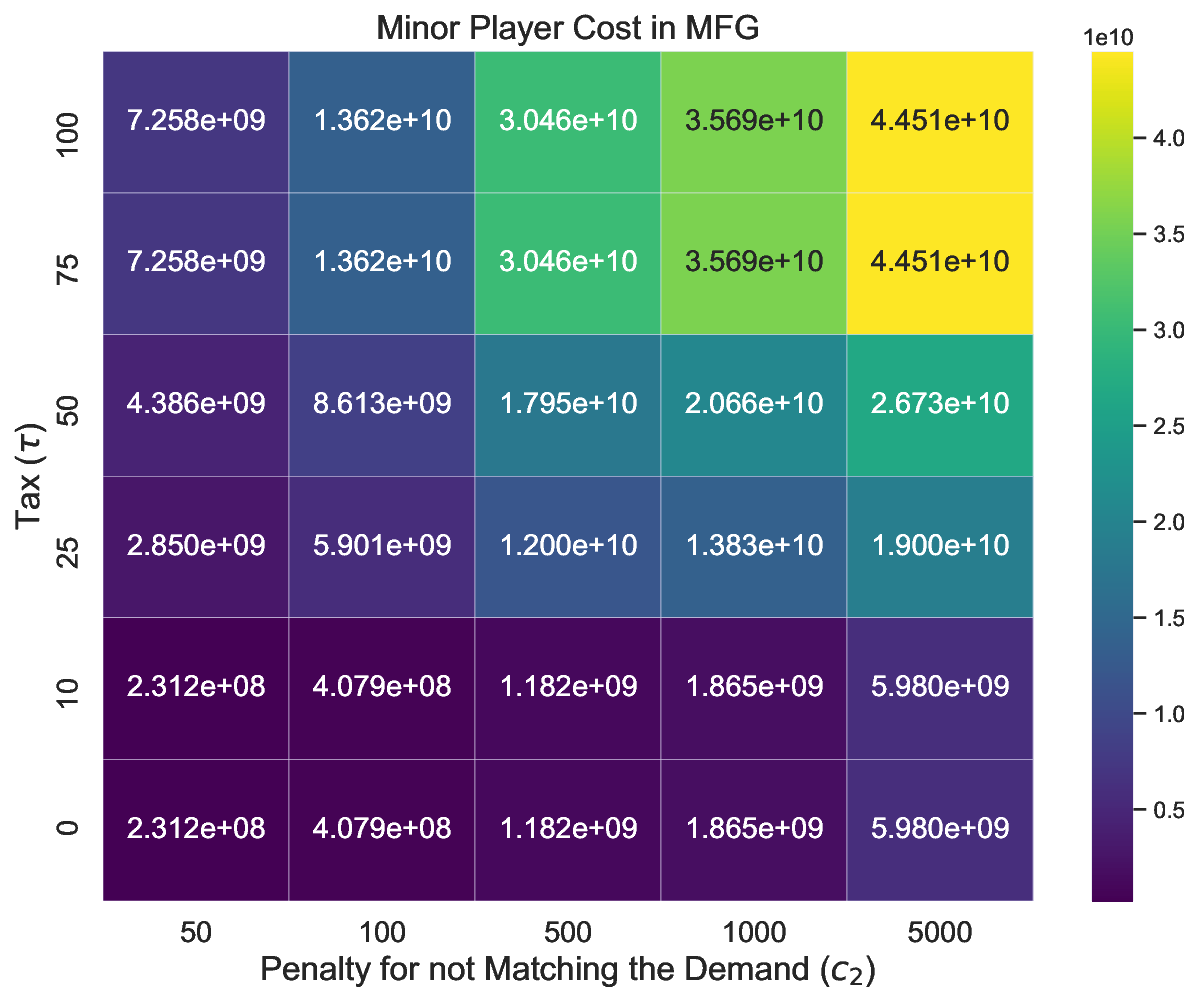}
\end{subfigure}
\begin{subfigure}{.44\textwidth}
    \includegraphics[width=1\linewidth]{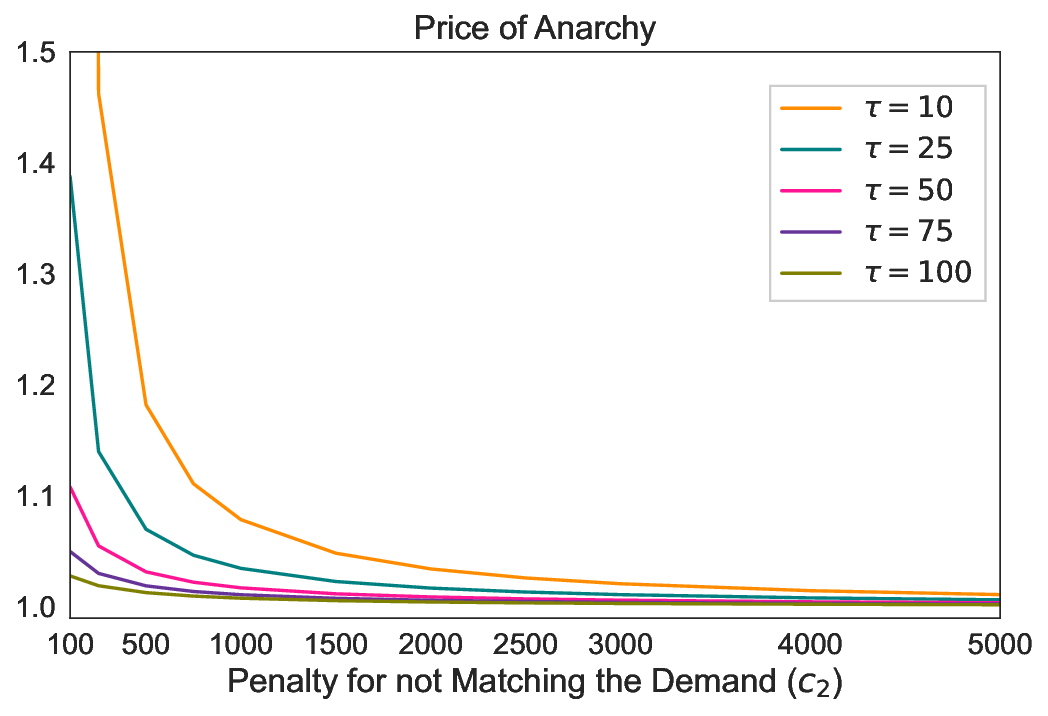}
\end{subfigure}
\caption{\textbf{Top:} Expected Cost of Minor Players in MFC (\textbf{left}), in MFG (\textbf{right}), \textbf{Bottom:} the Price of Anarchy given different penalty for not matching the demand and tax levels.}
\label{fig:minorcost_poa}
\end{figure}

From the heat maps in Figure \ref{fig:minorcost_poa}, we see that the expected cost of the representative producer is increasing with the carbon tax $\tau$ and the penalty $c_2$. The second observation is that as expected, for any given couple $(\tau, c_2)$, the expected cost is higher in the Nash equilibrium than for the Social Optimum. Next, we quantify how inefficient the Nash equilibrium is, and the effect of $\tau$ and $c_2$ on this inefficiency. In other words, we quantify the adverse effect of the non-cooperative behavior of the producers by computing the Price of Anarchy (PoA) defined in \eqref{eq:poa_minor} for different values of $\tau$ and $c_2$.
\begin{equation}\label{eq:poa_minor}
    PoA(\tau, c_2) = \frac{\inf_{N_t, R_e} C^{MFG}(N_t, R_e; \bar Q,\tau, c_2)}{\inf_{N_t, R_e} C^{MFC}(N_t, R_e; \bar Q, \tau, c_2)}.
\end{equation}
The results are given in the bottom subfigure in Figure \ref{fig:minorcost_poa}.
Since for any given $(\tau,c_2)$ the expected cost in a MFG equilibrium is higher, PoA is expected to be greater than $1$ and as it gets higher, the Nash Equilibrium is getting less efficient.
It can be seen that PoA gets smaller as we increase $\tau$ and $c_2$. This means that for higher levels of $\tau$ and $c_2$, the expected costs of to the producers become closer. In other words, \textit{the impact of the social planner diminishes and the advantages of cooperation lessen as the regulator imposes stricter regulations.}

\subsubsection{Electricity Production Decomposition Analysis}
\begin{figure}[H]
\centering
\begin{subfigure}{.99\textwidth}
    \includegraphics[width=1\linewidth]{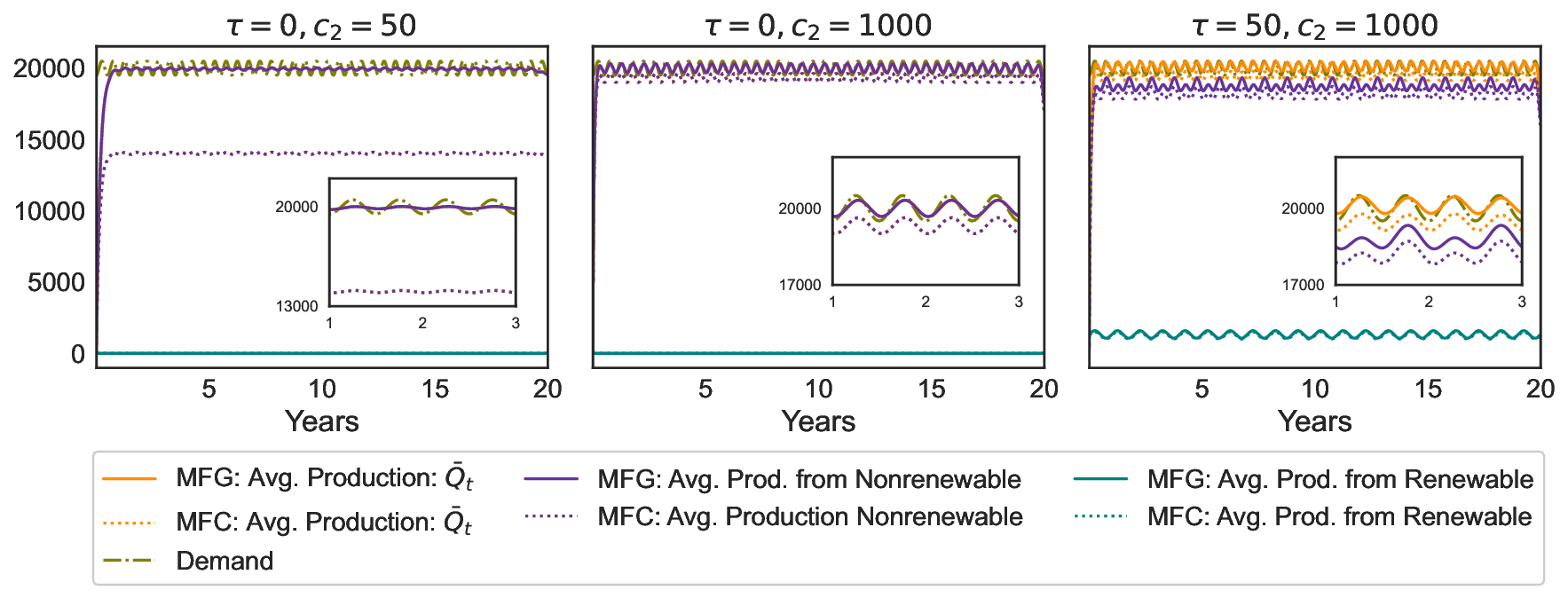}
\end{subfigure}
\caption{Production decompostion, \textbf{left:} when there is no tax ($\tau$) and the penalty for not matching the demand ($c_2$) is low, \textbf{middle:} when there is no tax and $c_2$ is high, \textbf{right:} when $\tau$ and $c_2$ are high. On the left and the middle subplots since all the production comes from the nonrenewable energy resources the average total production lines (colored in orange) are not seen.}
\label{fig:production}
\end{figure}

Here, we analyze the effect of the penalty $c_2$ for not matching the demand and the carbon tax $\tau$, on the optimal energy production portfolio in both MFC and MFG models. Figure \ref{fig:production} shows the total production and the decomposition of this production over a $20$ year period together with a detailed zoom in behavior between years $1$ and $3$.

The left subfigure in Figure~\ref{fig:production}, shows that the demand is not matched by the producers in the MFC case. This is because the penalty coefficient $c_2$  is low and the increased revenue from scarce supply is more advantageous. \textit{Here, we see that in the control setting, producers behave as a big monopoly when not matching the demand is inexpensive.} When the penalty is increased the middle subfigure in Figure~\ref{fig:production} shows that producers try to match the demand and their behaviors in the MFC and MFG cases are similar. In both of these figures there is no carbon tax, therefore the producers do not have incentives to invest in renewable energy, and as a result, all the production is exclusively from the nonrenewable sources. On the right subfigure in Figure~\ref{fig:production} \textit{when the carbon tax is increased we see that the producers have an incentive to invest in renewable energy.}

\begin{figure}[H]
\centering
\begin{subfigure}{.325\textwidth}
    \includegraphics[width=1\linewidth]{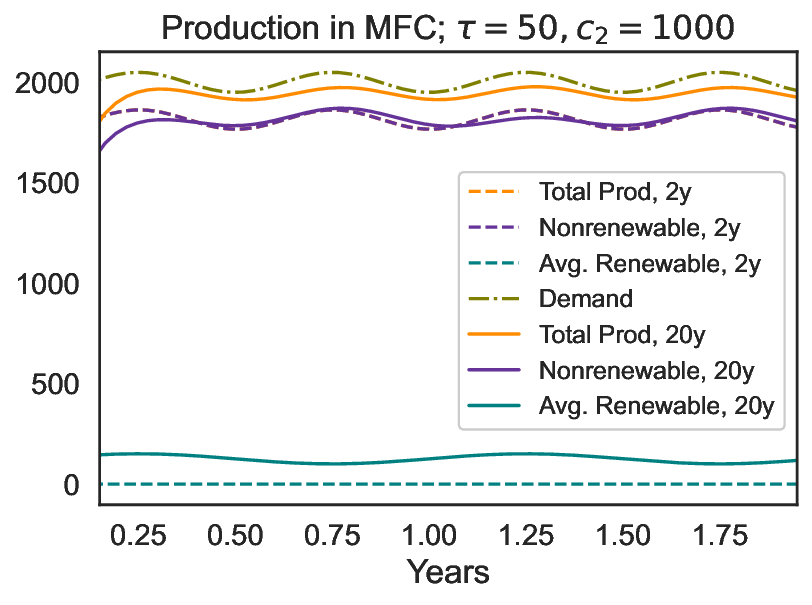}
\end{subfigure}
\begin{subfigure}{.325\textwidth}
    \includegraphics[width=1\linewidth]{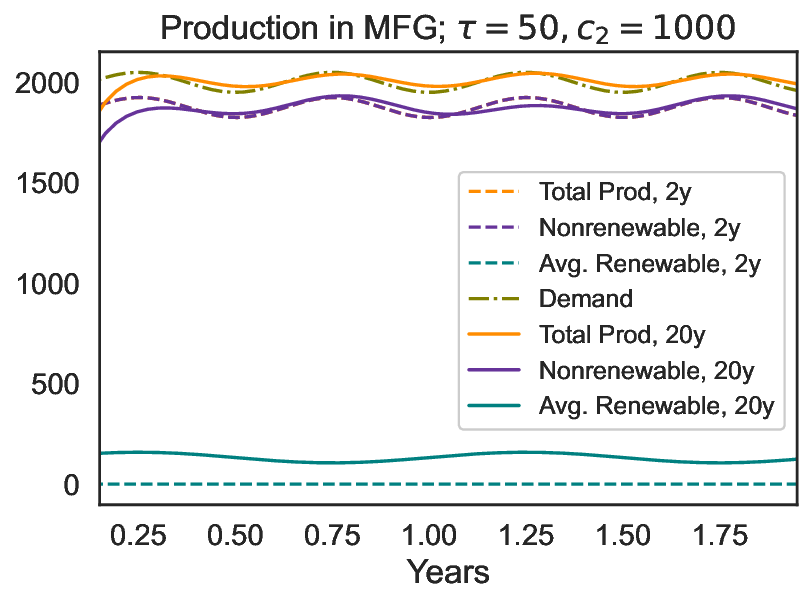}
\end{subfigure}
\begin{subfigure}{.325\textwidth}
    \includegraphics[width=1\linewidth]{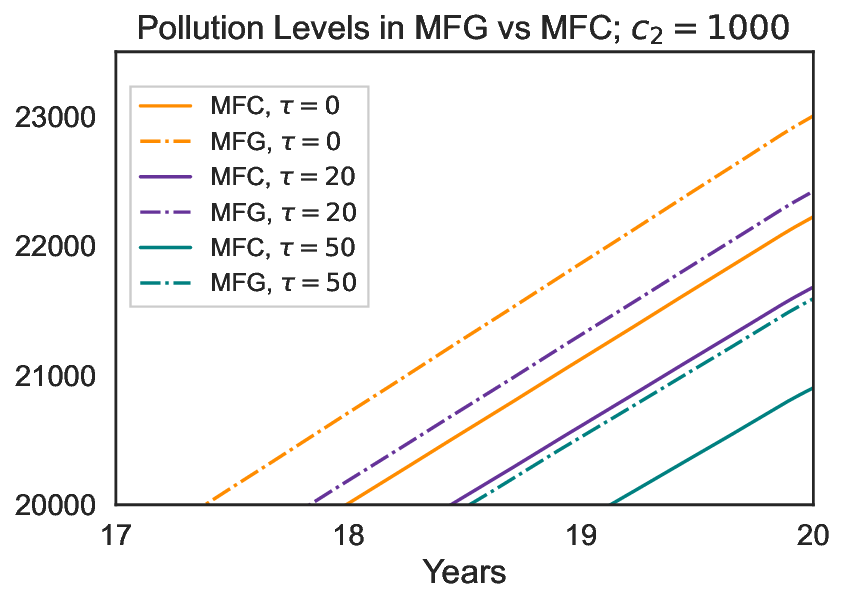}
\end{subfigure}
\caption{Planning time horizon effect in MFC (\textbf{left}), and in MFG (\textbf{middle}); End of the 20 years pollution in MFC and MFG (\textbf{right}).}
\label{fig:time_poll}
\end{figure}

We also analyze the effect of the planning horizon where we compare the cases in which the producers are planning for the next $2$ years vs. planning for the next $20$ years. As it can be seen in the left and middle subfigures in Figure~\ref{fig:time_poll}, when the planning horizon is short, the fixed costs of renewable energy outweigh its advantages. \textit{Short-sighted producers do not have an incentive to invest in renewable energy production.}

\subsubsection{Pollution Analysis}
The right subfigure in Figure~\ref{fig:time_poll}, shows that \textit{whatever the level of the carbon tax, the terminal pollution levels are higher when the producers are competitive} (MFG). Further,\textit{ in the absence of a carbon tax, producers can decrease the pollution levels further by cooperating and following a social planner instead of implementing a carbon tax.}

\section{Models with a Regulator}
\label{sec:regulator_model}
In this section, we describe how the previous models can be extended to include a major player in charge of choosing the tax level $\tau$ on behalf of a policy maker, and the penalty $c_2$ for not matching the demand on behalf of system operator. We shall treat this major player as a \textit{regulator}, and we shall often speak of minor players when we talk about the producers.
We extend the ``minor player only" model used previously by offering the producers the option to withdraw their entire production, de facto \textit{walking away} from the contract imposed by the regulator. This decision is made when the expected cost to the producer is higher than a fixed level above which producing at such a level of loss does not make sense. If we refer to the plots in Figure~\ref{fig:minorcost_poa}, we can see that the cost of the minor player is increasing with higher tax and the penalty for not matching the demand. Therefore, the regulator should be careful not to enact policies with very high values of $\tau$ and $c_2$. 

In the new model, the regulator does not have a private state per se. It only has $2$ controls which are the carbon tax level ($\tau \in \mathbb{R}_+$) and the penalty $(c_2\in \mathbb{R}_+)$.
Both controls are assumed to be time independent. This assumption is especially realistic when the period $[0,T]$ is too short for changes in regulation to make sense.
The\textit{ cost function} of the regulator is given as:
 \begin{equation}
 \begin{aligned}
 \label{eq:regulatorcost}
 J\Big(\tau, c_2; \bar P_T, \bar Q\Big)=&
 \underbracket{\alpha_1 \big( \bar P_T - \bar P^*_T \big)_+}_\text{Term~1} \underbracket{
 - \alpha_2 \tau \big(\bar P_T-\bar P_0\big)}_\text{Term~2} \underbracket{ +  \alpha_3 \big|\tau \big|^2}_\text{Term~3} \underbracket{ + \int_0^T  \alpha_4 \big|\bar{Q}_t - D_t \big|^2 dt}_\text{Term~4} \underbracket{+ \alpha_5c_2^2}_\text{Term~5}.
\end{aligned}
\end{equation}
The first term is the cost for exceeding the pollution target $\bar P_T^*\geq\bar P_0$ announced at $t=0$. Since we use the notation $x_+=\max(0,x)$, there is no penalty if the terminal pollution level is below the target. 
The constant $\alpha_1>0$ quantify the size of the penalty. 
The second term is the revenue from the carbon tax. 
To prevent the regulator from choosing an abusive high tax to increase their revenue, Term~3 is added to represent a reputation cost ($\alpha_2, \alpha_3>0$).
The joint roles of Term~4 and Term~5 is to insure that the responsibility of matching the demand is not only incumbent on the producers, but also on the regulator, influencing the choice of $\alpha_4>0$. This is consistent with our characterization of our major player / regulator as a policy maker as well as a system operator bearing the brunt of managing the ancillary services to avoid disruptions like \textit{system black-outs}.

\subsection{Equilibrium Notions}
\label{subsec:regulator_equilibrium}
We analyze two types of equilibria in the models with a regulator. In both cases, we consider that the regulator announces their policy first, and the producers react accordingly. This is in the realm of Stackelberg games. We call the first equilibrium \textit{Stackelberg MFC equilibrium}. In this case, the regulator assumes that a social planner chooses the controls used by the electricity producers. The latter behave like one big monopolistic firm. Therefore, the regulator chooses the tax level, $\tau$ and penalty coefficient $c_2$, assuming that the producers will settle in a MFC optimum. Note that in this interpretation the regulator and the social planner are two different entities. We define this equilibrium formally as: 

\begin{definition}[\textit{Stackelberg MFC equilibrium}]
    For every $(\tau, c_2)$, let\\
    $\Big(N^*(\tau, c_2), R^*_e(\tau,c_2)\Big)$  be the social planner's \textit{MFC optimum} given the tax level $\tau$ and the penalty coefficient $c_2$. In other words, for every $\tau, c_2$ and any admissible $\big(N, R_e\big)$, we have:
\begin{equation*}
\begin{aligned}
    &C\Big(\textcolor{Bittersweet}{\big(N, R_e\big)}; \bar{X}\big(\textcolor{Bittersweet}{N, R_e}\big), (\tau, c_2) \Big)
    \ge \\
    &\pushright{C\Big(\big(N^*(\tau, c_2), R^*_e(\tau, c_2)\big); \bar{X}\big(N^*(\tau, c_2), R^*_e(\tau, c_2)\big), (\tau, c_2)\Big)},
    \end{aligned}
\end{equation*}
where we added the notation $\bar{X}\big(N, R_e\big)$ to emphasize the parameters for which the mean field term $\bar X$ is computed.
Then the strategy profile $(\tau^*, c_2^*)$ is \textit{Stackelberg MFC equilibrium with a regulator} if, for any admissible $(\tau, c_2)$:
$$
    J\Big(\textcolor{Bittersweet}{(\tau, c_2)}; \bar{X}\big(N^*(\textcolor{Bittersweet}{\tau, c_2}),R^*_e(\textcolor{Bittersweet}{\tau, c_2}) \big)\Big)
    \ge
    J\Big((\tau^*, c_2^*); \bar{X}\big(N^*(\tau^*, c^*_4),R^*_e(\tau^*, c_2^*) \big)\Big).
$$
\end{definition}
The second equilibrium is called \textit{Stackelberg MFG Equilibrium}. In this one regulator assumes that electricity producers are competitive and it chooses $\tau$ and $c_2$ levels by assuming that the minor player population is at Nash Equilibrium. We can define this Equilibrium formally as 

\begin{definition}[\textit{Stackelberg MFG Equilibrium}]
    For every $(\tau, c_2)$, let \\
    $\Big(\hat N(\tau, c_2), \hat R_e(\tau,c_2)\Big)$  be the producers MFG \textit{Nash equilibrium} given the tax level $\tau$ and the demand satisfaction coefficient $c_2$. In other words, for any admissible $\big(N, R_e\big)$, we have:
\begin{equation*}
\begin{aligned}
    &C\Big(\textcolor{Bittersweet}{\big(N, R_e\big)}; {\bar{X}}\big(\hat N(\tau, c_2), \hat R_e(\tau, c_2)\big), (\tau, c_2) \Big)
    \ge \\
    &\pushright{C\Big(\big(\hat N(\tau, c_2), \hat R_e(\tau, c_2)\big); {\bar{X}}\big(\hat N(\tau, c_2), \hat R_e(\tau, c_2)\big), (\tau, c_2)\Big).}
\end{aligned}
\end{equation*}
Then the strategy profile $(\hat{\tau}, \hat{c}_2)$ is a \textit{Stackelberg MFG equilibrium with a regulator} if, for any admissible $(\tau, c_2)$, we have:
$$
    J\Big(\textcolor{Bittersweet}{(\tau, c_2)}; {\bar{X}}\big(\hat N(\textcolor{Bittersweet}{\tau, c_2}), \hat R_e(\textcolor{Bittersweet}{\tau, c_2}) \big)\Big)
    \ge
    J\Big((\hat{\tau}, \hat{c}_2); {\bar{X}}\big(\hat N(\hat{\tau},\hat{c}_2), \hat R_e(\hat{\tau},\hat{c}_2)\big)\Big).
$$
\end{definition}

\section{Numerical Results in the Presence of a Regulator}
\label{sec:reg_main_numeric_res}

\subsection{Algorithms}\label{subsec:stackelberg_alg}

To implement the walk-away option of the producers, we modify the \texttt{SocialOpt} and \texttt{NashEq} algorithms. This is done by simply adding an \texttt{IF} condition to these algorithms that assigns \texttt{Accept=1} if the cost of the minor player is lower than the threshold and \texttt{Accept=0} otherwise. After the algorithms for the producers are modified and called "\texttt{ModifiedSocialOpt}" and "\texttt{ModifiedNashEq}" respectively, we implement the Stackelberg Equilibrium Algorithm where we assume that if the producers reject the contract (\texttt{Accept=0}), the regulator cost is equal to infinity.

\begin{algorithm}[H]
\caption{\small Computation of Regulator's Cost} 
{\small\begin{algorithmic}[1]
\Function {\texttt{RegulatorCost}}{$\tau, c_2, \bar X, \text{Accept}$}
\vskip2mm
    \If{Accept = 1}
    \State Compute the regulator's cost $J=J(\tau, c_2; \bar X)$ by using \eqref{eq:regulatorcost}
    \vskip1mm
    \Else
    \State $J = \infty$ 
    \EndIf
    \vskip2mm
    \State \Return $J$
\EndFunction
\vskip3mm
\end{algorithmic}}
\end{algorithm}

\begin{algorithm}[H]
\caption{\small Search for a Stackelberg Equilibrium with MFC and MFG}  
{\small\begin{algorithmic}[1]
\Function {\texttt{StackelbergEq}}{Type} 
\vskip2mm
    \If{Type = MFC}
        \State Search for optimal $(\hat \tau, \hat c_2)$ couple where optimal mean field $(\tau, c_2) \rightarrow \bar X(\tau, c_2)$, investment in renewable $(\tau, c_2) \rightarrow R_e(\tau, c_2)$ and minor cost $(\tau, c_2) \rightarrow c(\tau,c_2)$ are computed by \texttt{ModifiedSocialOpt} algorithm and optimal cost of regulator $(\tau, c_2) \rightarrow J(\tau, c_2; \bar X(\tau, c_2), \text{Accept})$ is found by using \texttt{RegulatorCost}
        \vskip1mm
    \ElsIf{Type = MFG}
        \State Search for optimal $(\hat \tau, \hat c_2)$ couple where optimal mean field $(\tau, c_2) \rightarrow \bar X(\tau, c_2)$, investment in renewable $(\tau, c_2) \rightarrow R_e(\tau, c_2)$ and minor cost $(\tau, c_2) \rightarrow c(\tau,c_2)$ are computed by \texttt{ModifiedNashEq} algorithm and optimal cost of regulator $(\tau, c_2) \rightarrow J(\tau, c_2; \bar X(\tau, c_2), \text{Accept})$ is found by using \texttt{RegulatorCost}
    \EndIf
    \vskip2mm
    \State Let $(\hat \tau, \hat c_2) = \argmin_{\tau, c_2} J(\tau, c_2; \bar X(\tau, c_2))$, $\hat{\bar X }= \bar X(\hat \tau, \hat c_2)$, $\hat{R_e}= R_e(\hat \tau, \hat c_2)$, $\hat{c}= c(\hat \tau, \hat c_2)$ and $\hat J = J(\hat \tau,\hat c_2; \hat{\bar X})$
    \vskip2mm
    \State \Return $(\hat \tau, \hat c_2, \hat{\bar X}, \hat R_e, \hat c, \hat J)$
\EndFunction
\vskip3mm
\end{algorithmic}}
\end{algorithm}

\begin{remark}
   In the two Stackelberg equilibria, the numerical algorithms only differ in the solution of producers' problem. 
\end{remark}

\subsection{Numerical Experiments}

\subsubsection{Analysis of Regulator's Cost}
\begin{figure}[H]
\centering
    \begin{subfigure}{.325\textwidth}
        \includegraphics[width=1\linewidth]{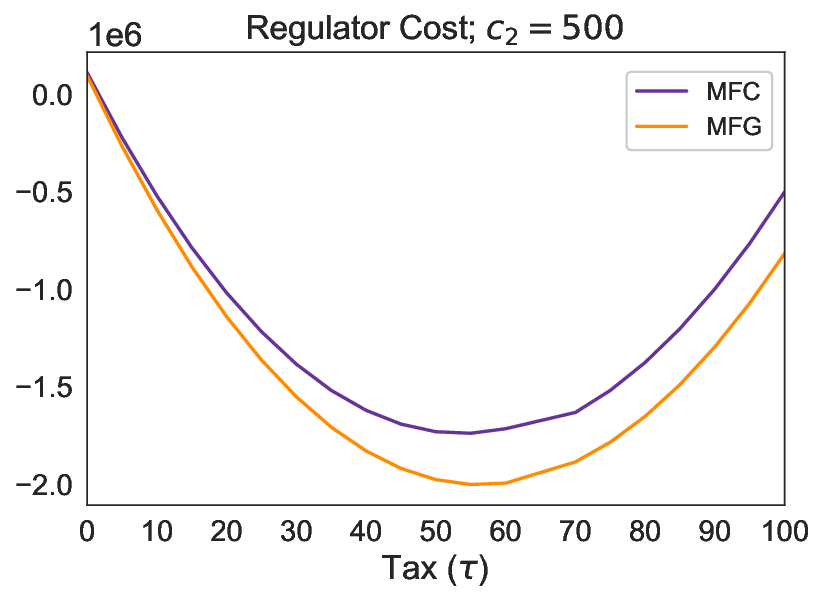}
    \end{subfigure}
    \begin{subfigure}{.325\textwidth}
        \includegraphics[width=1\linewidth]{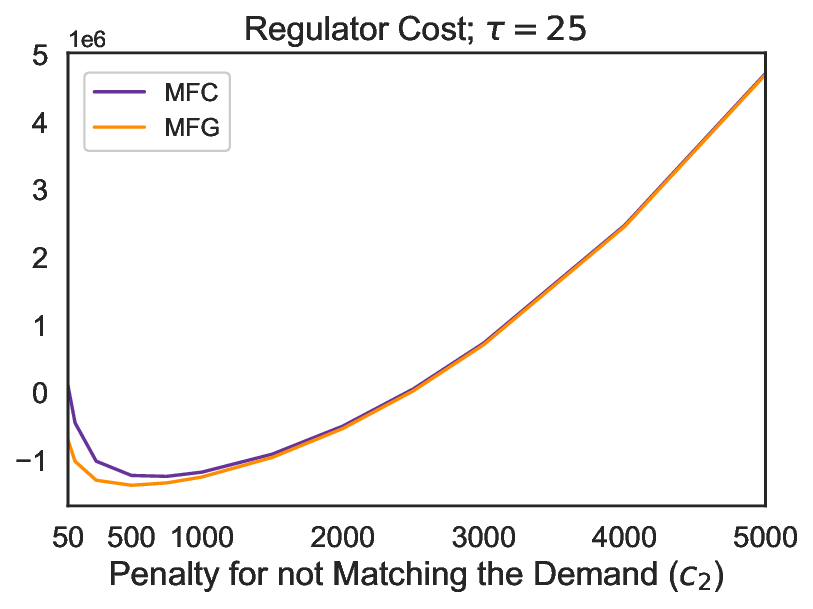}
    \end{subfigure}
    \caption{Regulator Cost in both MFC and MFG settings where penalty for not matching the demand, $c_2$, (\textbf{left}) or the carbon tax, $\tau$, (\textbf{right}) kept fixed.}
    \label{fig:reg_2d}
\end{figure}
\begin{figure}[H]
\centering
    \begin{subfigure}{.325\textwidth}
        \includegraphics[width=1\linewidth]{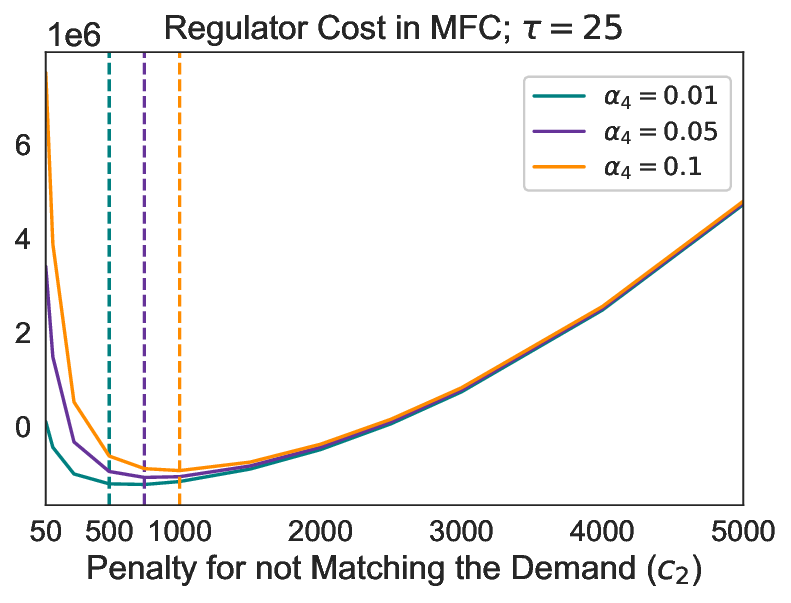}
    \end{subfigure}
    \begin{subfigure}{.325\textwidth}
        \includegraphics[width=1\linewidth]{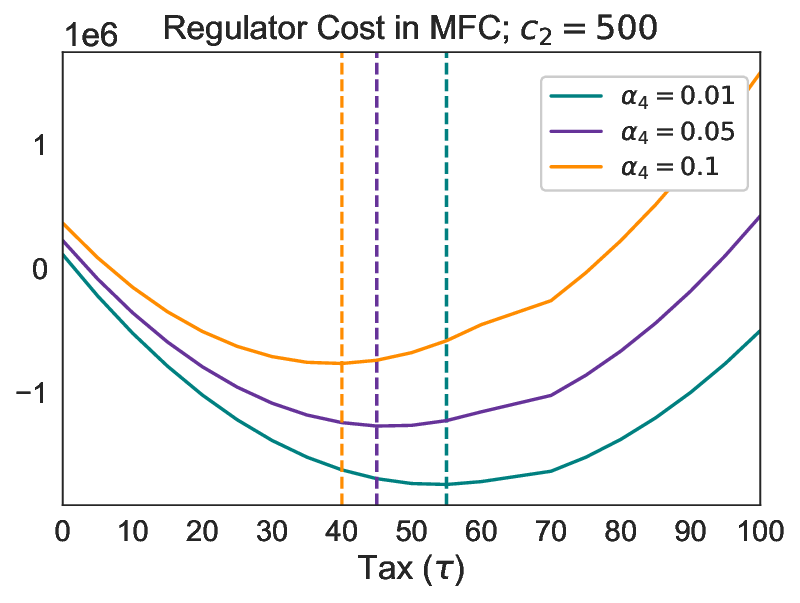}
    \end{subfigure}
    \begin{subfigure}{.325\textwidth}
        \includegraphics[width=1\linewidth]{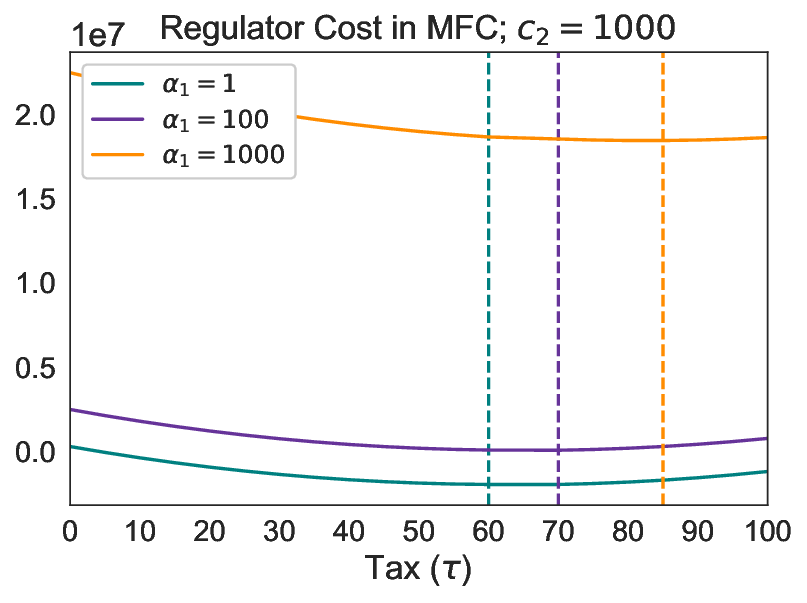}
    \end{subfigure}
    \caption{Effects of the coefficients of the regulator's cost and decisions: \textbf{left:} the effect of the importance given to the demand matching by the regulator, $\alpha_4$, on the decision of the penalty imposed on minor player for not matching the demand, \textbf{middle:} the effect of $\alpha_4$, on the decision of the carbon tax, \textbf{right:} the effect of the importance given to minimizing the excess pollution, $\alpha_1$, on the decision of the carbon tax. In the subplots the dashed vertical lines show the position of the minimizers of the regulator's cost with the corresponding color.}
    \label{fig:reg_2d_coeff}
\end{figure}

 
For the experiments of this section, we used the same parameters as for the producers' model in the previous section. For the regulator we used\footnote{For the minor player's problem, we have been able to choose realistic parameters by using real life data as explained in Subsection~\ref{subsec:minor_param}. However the parameters for the regulator's cost depend on the type of the regulator we focus on. For example, a regulator can care about minimizing the pollution relatively more than the other objectives or its main goal can be to maximize demand matching by the producers. Therefore, in the experiments, we focus on showing the effect of these different parameter choices on the decision of the regulator.}:
\vskip 12pt
{\small\begin{center}
{\renewcommand{\arraystretch}{1.2}
\begin{tabular}{c} 
 \hline
$\alpha_1 = 1$, $\alpha_2 = 3.3$, $\alpha_3 = 500$, $\alpha_4 = 0.01$, $\alpha_5 = 0.25$ \\
$\tau \in \{0, 10, 15, 20, 25, 30, 40, 50, 75, 100\}$ \\
$c_2 \in\{50, 100, 250, 500, 750, 1000, 1500, 2000, 2500, 3000, 4000, 5000\}$\\
 \hline
\end{tabular}}
\end{center}}
\vskip5mm
First, we analyze the regulator's expected cost for different values of the carbon tax given a fixed penalty for not matching the demand. Then we switch the roles of the two controls of the regulator. Plots in Figure~\ref{fig:reg_2d} show that \textit{the cost of the regulator is convex as a function of the carbon tax or the penalty, when the other control is fixed.} We also analyze the effect of the coefficients in the regulator's cost. First we start with the analysis of the importance given to demand matching by the regulator by tracking the effect of $\alpha_4$ in regulator's cost. The left subfigure in Figure~\ref{fig:reg_2d_coeff} shows that when the tax is fixed, the regulator's minimum cost is attained at higher $c_2$ values when $\alpha_4$ is higher. This shows that \textit{the regulator should impose higher penalties for not matching the demand to producers when demand matching is more important for the regulator.} The middle subfigure shows that when the penalty for not matching the demand is fixed and when $\alpha_4$ is higher then the optimal tax is lower. The reasoning here is that when the regulator cares about demand matching \textit{since the production from renewable energy is more unpredictable, the regulator is not opposed to the nonrenewable energy usage in order to have more stable demand matching.} Finally the right subfigure shows the effect of the importance given to minimizing the excess pollution by the regulator by tracking the effect of $\alpha_1$ in regulator's cost. Here, it can be seen that when the penalty for not matching the demand is fixed, the \textit{optimal carbon tax is higher when the regulator wants to keep pollution at a lower level.} These subfigures are for the MFC case but similar results hold in the MFG case as well.

\begin{figure}[H]
\centering
    \begin{subfigure}{.99\textwidth}
        \includegraphics[width=1\linewidth]{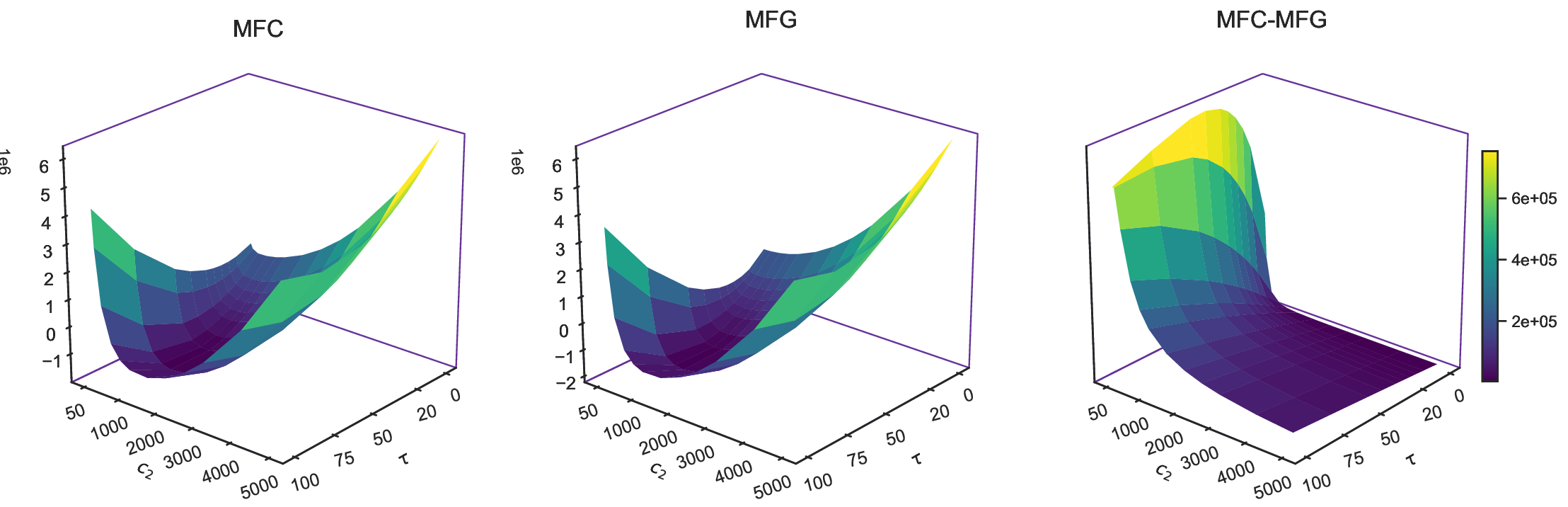}
    \end{subfigure}
    \caption{Regulator cost given admissible $c_2$ and $\tau$ values in when the regulator cares about matching the demand in MFC (in both MFC and MFG settings) and in MFG (\textbf{middle}); the difference of the regulator's cost between MFC and MFG given any admissible $c_2$ and $\tau$ couples (\textbf{right}).}
    \label{fig:reg_3d}
\end{figure}

Finally, Figure \ref{fig:reg_3d} gives 3D plots of the regulator cost as a function of their controls $\tau$ and $c_2$. Here, the minimum is attained at $(\tau^*, c_2^*) = (50, 1000)$ in the MFC case and at $(\hat \tau, \hat c_2) = (75, 1000)$ in the MFG case. Also, we see that for any given tax and penalty, the expected cost of the regulator is higher if the producers are cooperative instead of competitive when the regulator gives more importance to demand matching. This is because for any given couple $(\tau, c_2)$, in the cooperative setting producers are behaving like a big monopolistic firm and care less about matching the demand than in the competitive setting in order to maximize their revenues by keeping the prices higher. \textit{When demand matching is important for the regulator, the regulator benefits from the competition among the electricity producers even if this competition creates adverse effect for the producers themselves.}

\section{Conclusion}
\label{sec:conclusion}
In this paper, we investigate the behavior of rational electricity producers in the presence of a carbon tax. We analyze how they manage the trade-off between reliance on traditional and predictable fossil fuel power production assets which emit greenhouse gas and hence cost revenues because of the carbon tax, and the temptation to invest in clean energy production assets which will not be the source of emissions but which make matching the demand problematic because of the volatility of their output. We study a large population of producers in two different models, a first one in which they compete and hopefully reach a Nash equilibrium, and a second one in which they cooperate and rely on the solution of a centralized optimization problem. In a second set of models, we introduce a regulator choosing the level of the carbon tax in hope to control the overall emissions in the economy, and a penalty to be imposed on producers failing to meet the demand in hope to avoid power outages and reputation costs. Our models are based on recent progress in the theory and the numerical analysis of mean field games and mean field control problems.

We showed that when the producers cooperate, they are better off by behaving like a single monopolistic firm. However, if the regulator raises excessively the penalty to match the demand, they can take advantage of the competitive behavior of the producers. While our models remain stylized, they open the door to more complex models, e.g. involving time dependent policies the regulator could base on the response of the producers. Furthermore, our models could be the used to include more features of the energy markets such as storage, and the interactions between neighboring states or countries.

\appendix

\section{Proofs of Theorems} 

\begin{proof}[Proof of Theorem~\ref{theorem:fbsde_mfg}]
Assume that the strategy couple $(\hat N, \hat R_e)$ is optimal in the mean field game and the corresponding mean field flow is given as $\bar Q = \bar Q(\hat N, \hat R_e)$. Now, assume that the representative player deviates from the optimal strategy and uses $(\widehat{N}+ \epsilon \widecheck{N}, \widehat{R}_e+ \epsilon \widecheck{R}_e)$. Then:
\begin{equation}
    \begin{aligned}
    \label{eq:perturbedcostmfg}
        \dfrac{dC(\widehat{N}+ \epsilon \widecheck{N}, \widehat{R}_e+ \epsilon \widecheck{R}_e; \bar Q)}{d\epsilon}|_{\epsilon=0} =& \mathbb{E}\Big[\int_0^T \Big[2{c_{1} \widehat{N}_t\widecheck{N}_t} +
    {p_1 \widecheck{\tilde{N}}_t}+
    {2c_2(\widehat{Q}_t-D_t)\widecheck{Q}_t}\\ 
    &-
    {c_3\big(\rho_0 + \rho_1(D_t- \bar Q_t)\big)\widecheck{Q}_t}\Big] dt + 
    {2\tau \widehat{P}_T\widecheck{P}_T} + 
    {p^\prime (\widehat{R}_e)}\widecheck R_e\Big].
    \end{aligned}
\end{equation}
Furthermore we have the following dynamics:
\begin{equation*}
    \begin{alignedat}{2}
        d\widecheck{Q}_t &= \kappa_{1} \widecheck{N}_tdt +  \kappa_{2}\widecheck{R}_e\left( \alpha \cos(\alpha t) dt + (\theta-S_t) dt + \sigma_0 d\widecheck W_t\right)
        \qquad, &&d\widecheck{E}_t = \delta  \widecheck{N}_t dt,
        \\ 
        d\widecheck{P}_t &= \widecheck{E}_t dt,  &&d \widecheck{\tilde{N}}_t = \widecheck{N}_t dt,
    \end{alignedat}
\end{equation*}
with initial conditions: $\widecheck{Q}_0 = \widecheck{E}_0 = \widecheck{P}_0 = \widecheck{\tilde N}_0 = 0$.
We can introduce the adjoint variables with the following dynamics:
\begin{align*}
        dY^1_t &= \Big(-{2c_2(\widehat{Q}_t-D_t)} +
    {c_3\big(\rho_0 + \rho_1(D_t- \bar Q_t)\big)}\Big)dt + Z_t^{1,1}d\widetilde{W}_t + Z_t^{1,2}dW_t,    \quad &&Y_T^1 = 0 \\
        dY^2_t &= \Big(\kappa_2 \widehat{R}_e Y^1_t +Y^2_t\Big)dt + Z_t^{2,1}d\widetilde{W}_t + Z_t^{2,2}dW_t,    &&Y_T^2 = 0  \\
        dY^3_t &= -Y_t^4dt + Z_t^{3,1}d\widetilde{W}_t + Z_t^{3,2}dW_t, &&   Y_T^3 = 0   \\ 
        dY^4_t &= Z_t^{4,1}d\widetilde{W}_t + Z_t^{4,2}dW_t,    && Y_T^4 = 2\tau\widehat{P}_T \\
        dY^5_t &= -p_1 dt + Z_t^{5,1}d\widetilde{W}_t + Z_t^{5,2}dW_t, && Y_T^5 = 0.
\end{align*}        
Plugging these dynamics in the perturbed cost function \eqref{eq:perturbedcostmfg} and applying integration by parts:

\begin{equation*}
    \begin{aligned}
        &\dfrac{dC(\widehat{N}+ \epsilon \widecheck{N}, \widehat{R}_e+ \epsilon \widecheck{R}_e; \bar Q)}{d\epsilon}|_{\epsilon=0}\\ 
        &\qquad\qquad=\mathbb{E}\Big[\int_0^T \Big[2{c_{1} \widehat{N}_t\widecheck{N}_t}+
        {\big(-dY_t^5 + Z_t^{5,1}d\widetilde{W}_t + Z_t^{5,2}dW_t\big) \widecheck{\tilde{N}}_t}\\ 
        &\qquad\qquad\qquad\qquad +\big(-dY_t^1 +Z_t^{1,1}d\widetilde{W}_t + Z_t^{1,2}dW_t\big) 
        \widecheck{Q}_t\Big] dt + {Y^4_T\widecheck{P}_T} + 
        {p^\prime (\widehat{R}_e)}\widecheck R_e\Big]\\
        &\qquad\qquad=\mathbb{E}\Big[\int_0^T \widecheck{N}_t \Big(2{c_{1} \widehat{N}_t} +Y_t^5+Y_t^1\kappa_1+Y_t^3\delta\Big)dt\Big] \\
        &\qquad\qquad\qquad\qquad+ \mathbb{E}\Big[\int_0^T \left(\kappa_2  Y_t^1\Big(\alpha \cos(\alpha t) + (\theta {-S_t}) \Big) dt  + {p^\prime (\widehat{R}_e)} \right)\widecheck R_e\Big].
    \end{aligned}
\end{equation*}
By optimality, the above expression should be equal to 0 for any given $\widecheck N_t$ and $\widecheck R_e$; therefore:
\begin{equation}
\label{eq:optcond_fbsde}
        \widehat{N}_t = -\frac{Y_t^1\kappa_1+Y_t^3\delta+Y_t^5}{2c_1}\quad \text{and} \quad
        \hat R_e = (p^{\prime})^{-1} \Big(-\mathbb{E}\Big[\int_0^T \kappa_2Y_t^1\Big(\alpha \cos(\alpha t) + (\theta {-S_t}) \Big) dt \Big]\Big).
\end{equation}
For the proof of the FBSDE system in the mean field control setting, assume that strategy couple $(\widehat N, \widehat  R_e)$ is optimal in the mean field control problem and the corresponding mean field is given as $\widehat{\bar Q}$. Now assume that the representative player deviates from optimal strategy and uses $(\widehat N + \epsilon \widecheck N, \widehat R_e + \epsilon \widecheck R_e)$. Since in the mean field control version every player is going to deviate, the mean field term also changes to $\widehat{\bar Q} +\epsilon \widecheck{\bar Q}$. Then:
\begin{equation}
    \begin{aligned}
    \label{eq:perturbedcostmfc}
        \dfrac{dC({\widehat N}+ \epsilon \widecheck{N}, \widehat {R}_e+ \epsilon \widecheck{R}_e; \bar Q)}{d\epsilon}|_{\epsilon=0} &=\mathbb{E}\Big[\int_0^T \Big[2{c_{1} \widehat {N}_t\widecheck{N}_t} +
        {p_1 \widecheck{\tilde{N}}_t}+
        {2c_2(\widehat {Q}_t-D_t)\widecheck{Q}_t}\\
        &\pushright{-
        {c_3\big(\rho_0 + \rho_1(D_t- 2{\widehat{\bar Q}}_t)\big)\widecheck{Q}_t}\Big] dt + 
        {2\tau \widehat {P}_T\widecheck{P}_T} + 
        {p^\prime (\widehat {R}_e)}\widecheck R_e\Big].}
    \end{aligned}
\end{equation}
For the sufficiency, let us assume that $(N_t, R_e)$ is the Nash equilibrium. Then, we want to show by using FBSDE result, $\forall (\widecheck N_t, \widecheck R_e)$ we have:
\begin{equation*}
     C(N, R_e; \bar Q)-C(\check N, \check R_e; \bar Q)\leq 0.
\end{equation*}
Therefore we have:

\begin{align*}
    &C(N_t, R_e; \bar Q)-C(\check N_t, \check R_e; \bar Q)\\
    &\qquad= \mathbb{E}\Big[\tau(P_T^2 - \check P_T^2) + p(R_e)-p(\check R_e)\\
    &\qquad\qquad+\int_0^T \big(c_1(N_t^2 - \check N_t^2) + p_1(N_t -\check{\tilde N}_t)+ c_2 (Q_t^2-\widecheck Q_t^2)-2c_2(Q_t-\widecheck Q_t)D_t\\
    &\pushright{ -c_3\rho_0(Q_t-\widecheck Q_t)-c_3\rho_1(D_t- \bar Q_t)(Q_t-\widecheck Q_t)\big)dt\Big]}\\
    &\qquad\leq \mathbb{E}\Big[Y_T^4(P_T - \widecheck P_T) + p(R_e)-p(\check R_e)\\
    &\qquad\qquad+\int_0^T \big(c_1(N_t^2 - \check N_t^2) + p_1(N_t -\check{\tilde N}_t)+ c_2 (Q_t^2-\widecheck Q_t^2)-2c_2(Q_t-\widecheck Q_t)D_t\\
    &\pushright{ -c_3\rho_0(Q_t-\widecheck Q_t)-c_3\rho_1(D_t- \bar Q_t)(Q_t-\widecheck Q_t)\big)dt\Big]}\\
    &\qquad= \mathbb{E}\Big[Y_T^4(P_T - \widecheck P_T) + p(R_e)-p(\check R_e)\\
    &\qquad\qquad+\int_0^T -dY_t^1(Q_t-\widecheck Q_t) + \int_0^T [-2c_2Q_t(Q_t -\widecheck Q_t)+c_2(Q_t^2-\widecheck Q_t^2)+c_1(N_t^2-\widecheck N_t^2)\\
    &\qquad\qquad\qquad\qquad\quad{+p_1(\tilde N_t -\widecheck{\tilde N}_t)]dt} + \int_0^T Z_t^{1,1}(Q_t- \widecheck Q_t)d\widecheck W_t + \int_0^T Z_t^{1,2}(Q_t- \widecheck Q_t)d W_t\Big]\\
    &\qquad\leq \mathbb{E}\Big[Y_T^4(P_T - \widecheck P_T) + p(R_e)-p(\check R_e)\\
    &\qquad\qquad+\int_0^T -dY_t^1(Q_t-\widecheck Q_t) + \int_0^T [-2c_2Q_t(Q_t -\widecheck Q_t)+2c_2Q_t(Q_t -\widecheck Q_t)\\
    &\pushright{+2c_1N_t(N_t-\widecheck N_t)+p_1(\tilde N_t -\widecheck{\tilde N}_t)]dt\Big]}\\
    &\qquad=\mathbb{E}\Big[Y_T^4(P_T - \widecheck P_T) + p(R_e)-p(\check R_e)\\
    &\pushright{+\int_0^T -dY_t^1(Q_t-\widecheck Q_t)+\int_0^T 2c_1N_t(N_t-\widecheck N_t)dt+\int_0^T-dY_t^5(\tilde N_t -\widecheck{\tilde N}_t)\Big]}.
    \end{align*}
Now we apply integration by parts:
\begin{align*}
    &\mathbb{E}\Big[Y_T^4(P_T - \widecheck P_T) + p(R_e)-p(\check R_e)\\
    &\pushright{+\int_0^T -dY_t^1(Q_t-\widecheck Q_t)+\int_0^T 2c_1N_t(N_t-\widecheck N_t)dt+\int_0^T-dY_t^5(\tilde N_t -\widecheck{\tilde N}_t)\Big]}\\
    &=\mathbb{E}\Big[Y_T^4(P_T - \widecheck P_T) + p(R_e)-p(\check R_e) +\int_0^T (\kappa_1Y_t^1 + 2c_1 N_t +Y_t^5)(N_t-\widecheck N_t)dt \\
    &\pushright{ + \int_0^T \kappa_2 (R_e -\widecheck R_e)Y_t^1(\alpha\cos(\alpha t)+(\theta-S_t))dt}\Big]\\
    &=\mathbb{E}\Big[p(R_e)-p(\check R_e) +\int_0^T (\delta Y_t^3 + \kappa_1Y_t^1 + 2c_1 N_t +Y_t^5)(N_t-\widecheck N_t)dt \\
    &\pushright{ + \int_0^T \kappa_2 (R_e -\widecheck R_e)Y_t^1(\alpha\cos(\alpha t)+(\theta-S_t))dt}\Big].
\end{align*}    

By using the optimality conditions we have:
\begin{align*}
    &\delta Y_t^3 + \kappa_1Y_t^1 + 2c_1 N_t +Y_t^5=0,\\
    &p^{\prime}(R_e)=-\mathbb{E}\Big[\int_0^T \kappa_2 Y_t^1(\alpha \cos(\alpha t)+\theta-S_t)dt\Big].
\end{align*}

Therefore we have:
\begin{equation*}
    \begin{aligned}
    &\mathbb{E}\Big[p(R_e)-p(\check R_e) +\int_0^T (\delta Y_t^3 + \kappa_1Y_t^1 + 2c_1 N_t +Y_t^5)(N_t-\widecheck N_t)dt \\
    &\pushright{+ \int_0^T \kappa_2 (R_e -\widecheck R_e)Y_t^1(\alpha\cos(\alpha t)+(\theta-S_t))dt}\Big]\\
    &\qquad=p(R_e)-p(\widecheck R_e) -(R_e -\widecheck R_e)p^{\prime}(R_e)\\
    &\qquad \leq 0,
    \end{aligned}
\end{equation*}
by using the convexity of function $p(\cdot)$.
\end{proof}

\begin{proof}[Proof of Theorem~\ref{theorem:fbsde_existence_mfg}]
In order to show existence of the MFG equilibrium mean field flow system, we show the existence of the solution of the FBSDE given in \eqref{eq:fbsde}. We first fix an $R_e$, then solve the FBSDE and calculate a new $R_e$ by using the optimality condition in \eqref{eq:fbsde_opt_cond}. In other words, we need to show that there exists a fixed point for a function $f$, $f(R_e)= R_e$, by using Brouwer Fixed Point Theorem. Therefore, we need to show that $f:[0, R_e^{\text{max}}]\mapsto [0, R_e^{\text{max}}]$ is continuous in $R_e$. In order to simplify the notations we define $X_t := [Q_t, S_t, E_t, P_t, \tilde{N}_t]$, $Y_t := [Y^1_t, Y^2_t, Y^3_t, Y^4_t, Y^5_t]$ and 
$$
Z_t := 
\begin{bmatrix}
Z_t^{1,1}\quad & Z_t^{2,1}\quad & Z_t^{3,1}\quad & Z_t^{4,1}\quad& Z_t^{5,1}\\
Z_t^{1,2}\quad & Z_t^{2,2}\quad & Z_t^{3,2}\quad & Z_t^{4,2}\quad& Z_t^{5,2}\\
\end{bmatrix}^{\top}.
$$

Further, we define:
\begin{equation*}\arraycolsep3pt 
    K^x =  -\frac{1}{2c_1}\begin{bmatrix}
        \kappa^2_1 & 0 & \kappa_1\delta & 0 & \kappa_1\\
        0  &  0  &  0  &  0 &  0\\
        \kappa_1\delta  &  0  &  \delta^2  &  0 &  \delta\\
        0  &  0  &  0  &  0 &  0\\
        \kappa_1  &  0  &  \delta  &  0 &  1
    \end{bmatrix},
        L^x = -(K^y)^{\top} = \begin{bmatrix}
        0 &  -\kappa_2 R_e & 0 &0 & 0\\
        0 & -1 & 0 & 0 & 0\\
        0 & 0 & 0 & 0 & 0\\
        0 & 0 & 1 & 0 & 0\\
        0 & 0 & 0 & 0 & 0
        \end{bmatrix},
    O^y = \begin{bmatrix}
        -c_3 \rho_1& 0 & 0 & 0 & 0\\
        0 & 0 & 0 & 0 & 0\\
        0 & 0 & 0 & 0 & 0\\
        0 & 0 & 0 & 0 & 0\\
        0 & 0 & 0 & 0 & 0
        \end{bmatrix},
\end{equation*}
\begin{equation*}\arraycolsep3pt 
    L^y = \begin{bmatrix}
        -2c_2 &0 & 0 & 0 & 0\\
        0 & 0 & 0 & 0 & 0\\
        0 & 0 & 0 & 0 & 0\\
        0 & 0 & 0 & 0 & 0\\
        0 & 0 & 0 & 0 & 0
        \end{bmatrix},
    M_t^y =  \begin{bmatrix}
        (2c_2+c_3\rho_1)D_t\\
        0\\
        0\\
        0\\
        -p_1
    \end{bmatrix},  
    M^x_t =  \begin{bmatrix}
        \kappa_2 R_e\Big( \alpha cos(\alpha t) + \theta \Big)\\
        \theta\\
        0\\
        0\\
        0
    \end{bmatrix},
\end{equation*}
\begin{equation*}\arraycolsep3pt 
    S_T = \begin{bmatrix}
        0 & 0 & 0 & 0 & 0\\
        0 & 0 & 0 & 0 & 0\\
        0 & 0 & 0 & 0 & 0\\
        0 & 0 & 0 & 2\tau & 0\\
        0 & 0 & 0 & 0 & 0
        \end{bmatrix},
    \Sigma = \begin{bmatrix}
        \kappa_2 R_e \sigma_0 & 0\\
        \sigma_0 & 0\\
        0 & \sigma_1\\
        0 & 0\\
        0 & 0
    \end{bmatrix},
    \widetilde{W}_t = \begin{bmatrix}
        \widecheck W_t\\
        W_t
        \end{bmatrix}. 
\end{equation*}

When $R_e$ is fixed, we can write the FBSDE system as 
    \begin{equation}
        \begin{aligned}
        dX_t &= K^x Y_t + L^x X_t + M_t^x +\Sigma d\widetilde W_t, &X_0=x_0\\
        dY_t &= K^y Y_t + L^y X_t + O^y \bar X_t + M_t^y + Z_t d\widetilde W_t, \qquad &Y_T= S_T X_T.
        \end{aligned}
    \end{equation}
    
For the proof of continuity, we first focus on the mean processes: 
    \begin{equation}
        \begin{aligned}
        d\bar X_t &= K^x \bar Y_t + L^x \bar X_t + M_t^x,\\
        d\bar Y_t &= K^y \bar Y_t + (L^y + O^y) \bar X_t + M_t^y.
        \end{aligned}
    \end{equation}
    
By introducing ansatz $\bar Y_t = \bar A_t \bar X_t + \bar B_t$, we can decouple the forward backward ODE and end up with
    \begin{equation}
        d\bar X_t = (K^x \bar A_t + L^x) \bar X_t + (K^x\bar B_t+ M_t^x),
    \end{equation}
where $\bar A_t$ is the solution of the following matrix Riccati differential equation and $\bar B_t$ is the solution of the linear ODE system: 
    \begin{equation}
        \begin{aligned}
        \dot{\bar A}_t & = - \bar A_t K^x\bar A_t + K^y \bar A_t - \bar A_t L^x +L^y + O^y,\\
        \dot{\bar B}_t & = (K^y - \bar A_t K^x)\bar B_t + M_t^y - \bar A_t M_t^x.
        \end{aligned}
    \end{equation}

By using the general ODE continuity results with respect to parameters, we can analyze the continuity of $\bar A_t$ with respect to $R_e$. Since the solution of the matrix Riccati equation is bounded (conditions that give the boundedness of a matrix Riccati differential equation can be found in \cite{Jacobson_riccatiboundedness}), the necessary lipschitzness assumption holds and we can conclude that $\bar A_t$ is continuous in $R_e$. Further, again by using the general ODE continuity results and the fact that $\bar A_t$ is continuous in $R_e$, we can conclude that $\bar B_t$ is continuous in $R_e$. In this way, we have shown that $\bar X_t$ is continuous in $R_e$.

Now assume that $\bar X$ is exogenous and define $U_t:= \kappa_2 \big(\alpha cos(\alpha t) + \theta - S_t \big)$, $\Delta Y_t^1 := Y_t^{1, R_e}- Y_t^{1, \widetilde R_e}$. We can write:
    \begin{equation}
        \begin{aligned}
            \hat{R}_r^{R_e} - \hat{R}_r^{\widetilde R_e} &= (p^{\prime})^{-1} \Big(\mathbb{E}\Big[\int_0^T - \big(Y_t^{1, R_e}- Y_t^{1, \widetilde R_e}\big)U_tdt\Big]\Big)\\
            &\leq (p^{\prime})^{-1} \Big(\int_0^T  \mathbb{E}\Big[\big|\Delta Y_t^1 U_t\big|\Big]dt\Big)\\
            &\leq (p^{\prime})^{-1} \Big(\int_0^T  \Big(\mathbb{E}\big[|\Delta Y_t^1|^2\big]\Big)^{1/2} \Big(\mathbb{E}\big[|U_t|^2\big]\Big)^{1/2}dt\Big)\\
            &\leq (p^{\prime})^{-1} \Big(C_T \Big(\sup_{t \in [0, T]} \mathbb{E}\big[|\Delta Y_t^1|^2\big]\Big)^{1/2}\Big)\\
            &\leq (p^{\prime})^{-1} \Bigg(\tilde C_T\mathbb{E}\Bigg[\Big(\int_0^T \Big|(L^x- \widetilde{L}^x) X_t + (M^x- \widetilde{M}_t^x)+ (K^x- \widetilde{K}^y)Y_t \\&\pushright{+O^y (\bar X_t-\widetilde{\bar X}_t )\Big|dt\Big)^2 
            + \int_0^T |\Sigma - \widetilde \Sigma|^2 dt \Bigg]^{1/2}}\Bigg),
        \end{aligned}
    \end{equation} 
where the first inequality comes from the convexity of $p$, the second inequality comes from Cauch-Schwarz inequality and the last inequality is the result of \cite[Theorem 5.4]{hu2019wellposedness}. By using the continuity of $\bar X$ in $R_e$ and the continuity of $(p^{\prime})^{-1}(\cdot)$, as $R_e-\tilde R_e$ goes to 0 the upper bound goes to 0. Therefore, we can infer that $\hat{R}_r^{R_e} - \hat{R}_r^{\widetilde R_e}$ also goes to 0, which gives continuity.

Since, we assume that $(p^{\prime})^{-1}:\mathbb{R}\mapsto [0, R_e^{\max}]$, we also have $f: [0, R_e^{\max}] \mapsto [0, R_e^{\max}]$. We conclude the existence proof by using Brouwer Fixed Point Theorem.

Uniqueness can be concluded as follows:
Assume there exist two mean field game equilibria: $({N}, R_e, {\bar Q}) =({N}_t, R_e, \bar Q_t)_{t\in[0,T]}$ and $({N}^{\prime}, R_e^{\prime}, \bar Q^{\prime}) = ({N}^{\prime}_t, R_e^{\prime}, {\bar Q}_t^{\prime})_{t\in[0,T]}$ such that $\bar Q \neq \bar Q'$. Then the control processes $({N}, R_e)$ and $({N}^{\prime}, R_e^{\prime})$ should differ since if they are the same we would have the same state processes and the distributions would be the same. By using the definition of ``minimizer" of a cost functional, we have:
\begin{equation*}
    \begin{aligned}
    C(N, R_e; \bar Q)&\leq C(N^{\prime}, R^{\prime}_e; \bar Q), \qquad C( N^{\prime}, R^{\prime}_e; \bar Q^{\prime})&\leq C( N, R_e; \bar Q^{\prime}).
    \end{aligned}
\end{equation*}
By adding the two inequalities, we get:
\begin{equation}
\label{eq:uniqueness_ineq}
    \Big(C(N, R_e; \bar Q)- C( N, R_e; \bar Q^{\prime})\Big)-\Big(C(N^{\prime}, R^{\prime}_e; \bar Q)-C( N^{\prime}, R^{\prime}_e; \bar Q^{\prime})\Big)\leq 0.
\end{equation}
Now we use the fact that the drift and the volatility terms are independent of the state distribution, $\mathcal{L}(X)=\boldsymbol{ \mu}$. Therefore, in environment $\boldsymbol{ \mu}$
, the controlled path driven by $({N}^{\prime}, R^{\prime}_e)$ is ${\bar Q}^{\prime}$ and in environment 
$\boldsymbol{ \mu}^{\prime}$
, the controlled path driven by $( N, R_e)$ is ${\bar Q}$. By using this, we write:
\begin{equation*}
    \begin{aligned}
    &C(N, R_e; \bar Q)- C( N, R_e; \bar Q^{\prime})
     = c_3\rho_1\int_0^T \bar Q_t (\bar Q_t - \bar Q_t')dt.
    \end{aligned}
\end{equation*}
In the same way, we have:
\begin{equation*}
    \begin{aligned}
    C(N^{\prime}, R^{\prime}_e; \bar Q)-C( N^{\prime}, R^{\prime}_e; \bar Q^{\prime})=
    c_3\rho_1\int_0^T \bar Q_t' (\bar Q_t - \bar Q_t')dt.
    \end{aligned}
\end{equation*}
Therefore the expression on the left of the inequality \eqref{eq:uniqueness_ineq} becomes:
\begin{equation*}
    \begin{aligned}
     &  \Big(C(N, R_e; \bar Q)- C( N, R_e; \bar Q^{\prime})\Big)-\Big(C(N^{\prime}, R^{\prime}_e; \bar Q)-C( N^{\prime}, R^{\prime}_e; \bar Q^{\prime})\Big)
    = c_3\rho_1\int_0^T (\bar Q_t - \bar Q'_t)^2 dt
    >0.
    \end{aligned}
\end{equation*}
By contradiction, we conclude the uniqueness.
\end{proof}

\begin{proof}[Proof of Theorem~\ref{theorem:fbsde_mfc}]
We introduce the same adjoint variables as for the MFG FBSDE except that for $Y_t^1$ we take:
\begin{equation*}
    dY^1_t = \Big(-{2c_2(\widehat{Q}_t-D_t)} +
    {c_3\big(\rho_0 + \rho_1(D_t- 2\widehat{\bar Q}_t)\big)}\Big)dt + Z_t^{1,1}d\widecheck{W}_t + Z_t^{1,2}dW_t,   \qquad Y_T^1 = 0.
\end{equation*}
By plugging the adjoint variable dynamics in the perturbed cost \eqref{eq:perturbedcostmfc} and applying integration by parts, we end up with the same optimality conditions as in \eqref{eq:optcond_fbsde}. 
The sufficiency condition is proved by following the same ideas as in the proof of Theorem~\ref{theorem:fbsde_mfg}: we have
\begin{equation*}
    \begin{aligned}
    &C(N_t, R_e; \bar Q_t)-C(\check N_t, \check R_e; \check{\bar{Q}}_t)\\
    &\qquad=p(R_e)-p(\check R_e)-p^{\prime}(R_e)(R_e-\check R_e)\\
    &\pushright{+\mathbb{E}\Big[\int_0^T\big[c_3\rho_1(Q_t \bar Q_t- \check Q_t \check{\bar Q}_t)-2c_3\rho_1\bar Q_t(Q_t-\check Q_t)\big]dt\Big]}\\
    &\qquad=p(R_e)-p(\check R_e)-p^{\prime}(R_e)(R_e-\check R_e)-\int_0^T \big[c_3\rho_1(\bar Q_t+ \check{\bar Q}_t)^2 \big]dt\\
    &\qquad\leq 0.
    \end{aligned}
\end{equation*}    
which is obtained by using the convexity of the function $p(\cdot)$. 
\end{proof}

\begin{proof}[Proof of Theorem~\ref{theorem:fbsde_existence_mfc}]
The proof of the existence of the solution in the MFC case follows the same ideas of the proof of Theorem~\ref{theorem:fbsde_existence_mfg} and for the sake of space, it is omitted. To prove uniqueness of the MFC optimal mean field term, we introduce an auxiliary MFG which has the same FBSDE as the MFC problem and for which we prove uniqueness. To wit, we first focus on the mean field game problem that has the same dynamics as in \eqref{eq:minordynamics} and that has the following cost functional for an infinitesimal agent given a mean field flow $(\bar Q_t)_t$:
\begin{equation*}
    \begin{aligned}
    \mathbb{E}\Big[\int_0^T \Big[c_{1} |N_t|^2 +
    p_1 \tilde{N}_t+
    c_2|Q_t-D_t|^2 - 
    c_3\big(\rho_0 + \rho_1(D_t- 2\bar Q_t)\big)Q_t\Big] dt + 
    \tau|P_T|^2 + 
    p(R_e)\Big].
    \end{aligned}
\end{equation*}
Following the idea given in the proof of Theorem~\ref{theorem:fbsde_mfg}, the FBSDE system that characterizes the solution of this new game is found to be the same FBSDE that characterizes the solution of the mean field control. 

Uniqueness of the mean field flow of the new mean field game can be proved by using the approach given in the proof of Theorem~\ref{theorem:fbsde_existence_mfg} and it is omitted for the sake of space. This in turn concludes the uniqueness for the mean field control problem.
\end{proof}

\begin{proof}[Proof of Lemma~\ref{lem:ode_mfg}]
Following \cite[ch. 3]{carmona_2018}, we write the Hamiltonian:
\begin{equation}
\label{Hamiltonian}
H(t, N, X, \bar{X}, q) = (Ax)^{\top} q + (B \cdot N) ^{\top} q + C_t^{\top} q + \frac{R}{2} |N|^2 +H_t^{\top} X + \bar{X}^{\top} F X + X^{\top} G X +J_t,
\end{equation}
where $q$ is the adjoint process. Therefore, the optimal $N$ to optimize $H$ can be expressed as:
\begin{equation}
\label{MinorOptimalControl}
    \hat{N}(q) = - R^{-1} B^{\top} q.
\end{equation}
By plugging $\hat N(q)$ in the Hamiltonian \eqref{Hamiltonian}, the \textit{optimal Hamiltonian} can be written as:
\begin{equation}
\label{eq:OptimalHamiltonian}
    \hat H(t, X, \bar{X}, q) = -\frac{1}{2} q^{\top} B R^{-1}B^{\top} q + X^{\top} A^{\top} q + C_t^{\top} q + H_t^{\top} X + \bar{X}^{\top} F X+ X^{\top} G X +J_t.
\end{equation}
Then the Hamilton Jacobi Bellman (HJB) equation can be written as:
\begin{equation}
\label{eq:hjb_mfg}
    -\frac{\partial u(t,X)}{\partial t} - tr(aD^2u(t,X))  = \hat H(t, X, \bar{X}_t, Du(t,X)),\qquad  u(T,X) = X^{\top} S_T X + p_2 R_e,
\end{equation}
where $\bar{X}_t = \int_{\mathbb{R}^5} X m(t,X)dX$. 
The Kolmogorov Fokker Planck (KFP) equation for our MFG problem can be written as:
\begin{equation}
\begin{aligned}
    \label{eq:kfp_mfg}
    &\dfrac{\partial m(t,X)}{\partial t} - tr(aD^2 m(t,X))+ \nabla_x\big(m(t,X)(AX-Br^{-1}B^{\top}Du(t,X)+C_t)\big)=0, \\
    &\pushright{\quad m(0,X) = m_0(X).}
\end{aligned}
\end{equation}
Introduce the following ansatz for the value function:
\begin{equation}
    \label{eq:value_ansatz}
     u(t,X) = \frac{1}{2} X^{\top} {\eta_t} X + X^{\top} {r_t} +s_t.
\end{equation}
We have $Du(t,X) = \eta_t X + r_t$ and $D^2u(t,x)=\eta_t$.
By plugging \eqref{eq:value_ansatz} into the HJB equation in \eqref{eq:hjb_mfg}, we obtain that $\eta_t$ is the solution of the following symmetric matrix Riccati equation:
\begin{equation}
    \label{eq:riccati_mfg}
    \frac{d{\eta_t}}{dt} - {\eta_t} BR^{-1}B^{\top} {\eta_t} + A^{\top} {\eta_t} + {\eta_t} A +2G = 0, \qquad\qquad {\eta_T} = 2S_T.
\end{equation}
This Riccati equation has a unique positive symmetric solution, see~\cite[ch. 14.3]{kucera_2009}. 
By plugging \eqref{eq:value_ansatz} into the HJB equation~\eqref{eq:hjb_mfg}, we also obtain the differential equation for $r_t$ that is coupled with $\bar X_t$:
\begin{equation}
    \label{eq:ode_r_mfg}
    -\frac{d{r_t}}{dt} = \left(A^{\top} - {\eta_t} B R^{-1} B^{\top}\right) {r_t} + {\eta_t} C_t + H_t + F^{\top}{\bar{X}_t}, \qquad\qquad {r_T} = 0.
\end{equation}
The differential equation for $\bar X_t$ can be found by plugging the ansatz in the KFP equation:
\begin{equation}
    \begin{aligned}
    \label{eq:KFP_mfg_ansatz}
        \frac{d \bar X_t}{dt}&=\frac{d}{dt}\int_{\mathbb{R}^5} X \cdot m(t,X) dX= \int_{\mathbb{R}^5} X \cdot \frac{\partial m(t,X)}{\partial t} dX\\
        &= \int_{\mathbb{R}^5} X \cdot \Big(tr(aD^2 m(t,X))- \nabla_x\big(m(t,X)(AX-BR^{-1}B^{\top}(\eta_tX +r_t)+C_t)\big)\Big) dX\\
        &=-\int_{\mathbb{R}^5} X\cdot \frac{\partial m(t,X)}{\partial X} (AX-BR^{-1}B^{\top}(\eta_tX +r_t)+C_t) dX \\
        &\pushright{- \int_{\mathbb{R}^5} X\cdot m(t,X) (A-BR^{-1}B^{\top}\eta_t)dX}\\
        &=(A-BR^{-1}B^{\top}\eta_t) \bar X_t -BR^{-1}B^{\top}r_t+C_t.
    \end{aligned}
\end{equation}
Finally from the HJB equation where the ansatz is plugged in we find that:
\begin{equation*}
    \frac{ds_t}{dt}= -tr(a\eta_t) + \frac{1}{2}{r_t}^{\top} B R^{-1} B^{\top} {r_t} -C_t^{\top} {r_t}- J_t, \qquad\qquad s_T=p_2R_e.
\end{equation*}
Therefore we have:
\begin{equation}
    \label{eq:s_mfg}
    s_t = p_2 R_e + \int_t^{T} \Big(tr(a{\eta_s}) -\frac{1}{2} {r_s}^{T} B R^{-1} B^{\top} {r_s} +C_s^{\top} {r_s}+ J_s\Big)ds.
\end{equation}
The expected cost of the representative minor player given fixed mean field and $R_e$ can be calculated by using:
\begin{align*}
\begin{split}
    \inf_{(N_t)_t} \tilde{C}^{MFG}\Big(N; R_e, \bar X \Big) &= \mathbb{E}\left[u(0, X_0)\right]\\
    &= \mathbb{E}\left[ \frac{1}{2} X_0^{\top}\eta_0 X_0 + X_0^{\top} r_0 +s_0\right]\\
    &= \frac{1}{2} \left( Var(\sqrt{\eta_0} X_0) + \mathbb{E}[\sqrt{\eta_0} X_0]^2 \right) +\bar X_0^{\top} r_0 + s_0.
\end{split}    
\end{align*}
\end{proof}

\begin{proof}[Proof of Theorem~\ref{the:exist_uniq_mfg}]
For the existence and uniqueness proof, we make use of Banach Fixed Point theorem. We follow the line of proof used for a stochastic system in \cite[Thm 5.1]{ma_2007}. First we fix $(r_t^1)_t$ and $(r_t^2)_t$, then corresponding $(\bar X^i_t)_t$ can be found by solving the following ODE:
\begin{equation*}
    d\bar X_t^i = [(A-BR^{-1}B^{\top}\eta_t) \bar X_t^i -BR^{-1}B^{\top} r_t^i + C_t] dt, \qquad \bar X_0^i=\bar x_0, \qquad i=\{1,2\}.
\end{equation*}
Further let $\widetilde{\bar X}_t = \bar X_t^1 - \bar X_t^2$ and $\tilde r_t = r_t^1 - r_t^2$, then we have:
\begin{equation*}
    d\widetilde{\bar X}_t = [(A-BR^{-1}B^{\top}\eta_t) \widetilde{\bar X}_t -BR^{-1}B^{\top} \tilde r_t] dt, \qquad \widetilde{\bar X}_0=0.
\end{equation*}

Now we introduce $(r_t^{i'})_t$ that solves:
    \begin{equation*}
            dr^{i'}_t = [(\eta_t BR^{-1}B^{\top}-A^{\top})r^i_t -\eta_t C_t -H_t -F^{\top}\bar X_t^i]dt, \qquad r^{i'}_T = 0,  \qquad\forall i \in \{1,2\} .
    \end{equation*}
    and let $\tilde r_t^{'}=r_t^{1'}-r_t^{2'}$. Then we have the following ODE:
    \begin{equation*}
            d \tilde r^{'}_t = [(\eta_t BR^{-1}B^{\top}-A^{\top})\tilde r_t -F^{\top}\widetilde{\bar X}_t]dt, \qquad \tilde r^{'}_T=0.
    \end{equation*}
    Therefore we have defined a mapping $r \mapsto r^{'}$. We now show that it is a contraction mapping to be able to use the Banach Fixed Point theorem. 
    
    \noindent {\bf Step 1. } Using It\^o's formula, we write the dynamics for $||\widetilde{\bar X}_t||^2$:
\begin{equation*}
    d||\widetilde{\bar X}_t||^2 = 2(\widetilde{\bar X}_t)^{\top} d \widetilde{\bar X}_t = 2(\widetilde{\bar X}_t)^{\top} [(A-BR^{-1}B^{\top}\eta_t) \widetilde{\bar X}_t -BR^{-1}B^{\top} \tilde r_t]dt.
\end{equation*}
By using these dynamics we can find a bound for $||\widetilde{\bar X}_t||^2$:
\begin{align}
    \label{eq:mfg_ode_xbar_bound}
    ||\widetilde{\bar X}_t||^2 &= \int_0^t 2(\widetilde{\bar X}_s)^{\top} [(A-BR^{-1}B^{\top}\eta_s) \widetilde{\bar X}_s -BR^{-1}B^{\top} \tilde r_s]ds\nonumber \\
   & \leq \int_0^t  (|| 2(A-BR^{-1}B^{\top}\eta_s) ||) ||\widetilde{\bar X}_s||^2 ds + \int_0^t  || BR^{-1}B^{\top}|| 2 <\widetilde{\bar X}_s, \tilde r_s> ds\nonumber\\
   & \leq \int_0^t  (|| 2(A-BR^{-1}B^{\top}\eta_s) ||) ||\widetilde{\bar X}_s||^2ds + \int_0^t  || BR^{-1}B^{\top}||\big(||
       \widetilde{\bar X}_s||^2 + ||\tilde r_s||^2\big)ds\nonumber\\
   & \leq \exp\Big(\int_0^t(2||A-BR^{-1}B^{\top}\eta_s|| +||BR^{-1}B^{\top}||)ds\Big) \int_0^t (||BR^{-1}B^{\top}||) ||\tilde r_s||^2 ds\nonumber\\
   & \leq C^{(1)} \int_0^T  ||\tilde r_s||^2 ds,
\end{align} 
    where the third to last inequalities stem from the Gronwall's inequality, and we define $C^{(1)} = \exp\Big(T\big(2||A||+ 2(||BR^{-1}B^{\top}||)||\eta||_T +||BR^{-1}B^{\top}||\big)\Big) ||BR^{-1}B^{\top}|| $ with $||\eta||_T := \sup_{0\leq t \leq T} ||\eta_t||$.

    \noindent {\bf Step 2. } We write the dynamics for $||\tilde r_t^{'}||^2$:
    \begin{equation*}
        d||\tilde r^{'}_t||^2 = 2(\tilde r^{'}_t)^{\top} d \tilde r^{'}_t = 2(\tilde r^{'}_t)^{\top} [(\eta_t BR^{-1}B^{\top}-A^{\top})\tilde r_t -F^{\top}\widetilde{\bar X}_t]dt.
    \end{equation*} 
    Now we find a bound for $||\tilde r_t^{'}||^2$ as follows by using Young's inequality:
    \begin{equation*}
        \begin{aligned}
           ||\tilde r^{'}_t||^2&= 
        \int_t^T 2(\tilde r^{'}_s)^{\top} [(A^{\top}-\eta_sBR^{-1}B^{\top}) \tilde r_s +F^{\top} \widetilde{\bar X}_s]ds\\
           &  \leq \int_t^T  ||A^{\top}-\eta_sBR^{-1}B^{\top} || 2<\tilde r^{'}_s, \tilde r_s> ds + \int_t^T  ||F^{\top}|| 2 <\tilde r^{'}_s, \widetilde{\bar X}_s> ds\\
           &  = \int_t^T  (||A^{\top}-\eta_sBR^{-1}B^{\top} || + ||F^{\top}|| ) ||\tilde r^{'}_s||^2 ds\\
           &\pushright{+ \int_t^T (||A^{\top}-\eta_sBR^{-1}B^{\top} ||)||\tilde r_s||^2 ds + \int_t^T  (||F^{\top}||) ||\widetilde{\bar X}_s||^2 ds.}
        \end{aligned}
    \end{equation*} 
    Now the expression found in \eqref{eq:mfg_ode_xbar_bound} can be plugged in and by using Gronwall's inequality:
    \begin{align*}
           ||\tilde r^{'}_t||^2 &\leq \int_t^T  (||A^{\top}-\eta_sBR^{-1}B^{\top} || + |F^{\top}|| ) ||\tilde r^{'}_s||^2 ds + \int_0^T (||A^{\top}-\eta_sBR^{-1}B^{\top} ||)||\tilde r_s||^2 ds +\\
            & \pushright{ \int_0^T  ||F^{\top}||  C^{(1)}  \Big(\int_0^T  ||\tilde r_s||^2 ds\Big) ds}\\
            &  \leq \int_t^T  (||A^{\top}-\eta_sBR^{-1}B^{\top} || + |F^{\top}|| ) ||\tilde r^{'}_s||^2 ds + \int_0^T (||A^{\top}-\eta_sBR^{-1}B^{\top} ||)||\tilde r_s||^2 ds +\\
            &  \pushright{ T\Big[ ||F^{\top}||  C^{(1)}  \Big(\int_0^T  ||\tilde r_s||^2 ds\Big) \Big]}\\
            &  \leq \exp\Big({T(||A|| +(||BR^{-1}B^\top||)||\eta||_T +||F^\top||)}\Big)\\
            &  \pushright{\int_0^T \Big[||A^{\top}-\eta_sBR^{-1}B^{\top} || + T||F^{\top}||  C^{(1)}  \Big]||\tilde r_s||^2 ds}.
    \end{align*} 
    
    Now define $||r||_T := \sup_{0\leq t \leq T} ||r_t||$, then we have:
       $ 
       ||\tilde r^{'}||_T^2 \leq c_T||\tilde r||_T^2
       $, 
    where $c_T$ is 
\begin{equation*}
\begin{aligned}
    c_T = &T e^{{T(||A|| +(||BR^{-1}B^\top||)||\eta||_T +||F^\top||)}}\\
    &\pushright{\times \big(||A^\top|| +(||\eta||_T +   T||F^{\top}||  e^{T(2||A|| +(2||\eta||_T+1 )||BR^{-1}B^{\top}||)})||BR^{-1}B^\top||\big).}
\end{aligned}
\end{equation*}
With small $T$ we have $c_T<1$, which concludes the proof.
\end{proof}

\begin{proof}[Proof of Lemma~\ref{lem:ode_mfc}]
For Mean Field Control problems we have the following HJB and FP systems. The detailed proof of the derivation can be found in \cite[ch. 6]{frehse_2013}.
\begin{equation*}
    \begin{aligned}
    &-\frac{\partial u(t,X)}{\partial t} - tr(aD^2u(t,X))  = \hat H(t, X, m(t), Du(t,X))\hfill\\
    &\pushright{+ \int_{\mathbb{R}^n} \frac{\partial \hat H}{\partial m}(t,\xi,m(t), Du(t,\xi))(x)m(t,\xi)d\xi,}\\
    &\pushright{u(T,X) = X^{\top} S_T X + p_2 R_e} \\
    &\dfrac{\partial m(t,X)}{\partial t} - tr(aD^2 m(t,X))+ \nabla_x\big(m(t,X)(AX-BR^{-1}B^{\top}Du(t,X)+C_t)\big)=0,\\ 
    &\pushright{m(0,X)=x_0,}
    \end{aligned}
\end{equation*}
where $\frac{\partial \hat H}{\partial m}$ denotes the G\^ateaux differential of $\hat H$ on $L^2(\mathbb{R}^5)$. 
As it can be seen, the KFP equation stays same as in MFC but HJB equation changes. By rewriting \eqref{eq:OptimalHamiltonian} as:
\begin{equation}
    \hat H(t, X, m, q) = -\frac{1}{2} q^{\top} B R^{-1}B^{\top} q + X^{\top} A^{\top} q + C_t^{\top} q + H_t^{\top} X + \Big(\int_{\mathbb{R}^n} \xi m(\xi)d\xi\Big)^{\top} F X+ X^{\top} G X +J_t.
\end{equation}
We find that
\begin{equation}
\label{eq:exp_hamiltonian_dens_der}
    \begin{aligned}
    \int_{\mathbb{R}^n}\frac{\partial \hat H(t, X, m, q)}{\partial m}(x)m(\xi)d\xi  = \int_{\mathbb{R}^n}X^{\top} F \xi m(\xi)d\xi=X^{\top} F \bar X .
    \end{aligned}
\end{equation}
Therefore, the HJB equation becomes:
\begin{equation}
\label{eq:hjb_mfc}
    -\frac{\partial u}{\partial t} - tr(aD^2u)  = \hat H(t, X, \bar{X}, Du) + X^{\top} F \bar X, \qquad  u(X_T,T) = X_T^{\top} S_T X_T + p_2 R_e .
\end{equation}
We introduce the same ansatz as in \eqref{eq:value_ansatz}. By plugging this ansatz in the HJB equation given in \eqref{eq:hjb_mfc}, we end up with the same Riccati equation and equation for $s_0$. Only the differential equation of $r_t$ changes as follows:
\begin{equation*}
    -\frac{d{r_t}}{dt} = \left(A^{\top} - {\eta_t} B R^{-1} B^{\top}\right) {r_t} + {\eta_t} C_t + H_t + F^{\top} {\bar{X}_t} + \textcolor{Bittersweet}{F {\bar{X}_t}}, \qquad {r_T} = 0, 
\end{equation*}
where $\bar{X}_t = \int_{\mathbb{R}^5} X m(t,X) dX$. 
Since the Kolmogorov-Fokker-Planck equation stayed the same, we have the same expression for the differential equation of $\bar X_t$ as in \eqref{eq:KFP_mfg_ansatz}. 
We obtain the MFC cost given any fixed $R_e$ by using the ansatz (see e.g.~\cite{AMSnotesLauriere} for more details):
\begin{align*}
\begin{split}
    \inf_{(N_t)_t} \tilde{C}^{MFC}\Big(N; R_e \Big) &= \mathbb{E}\left[u(0, X_0)\textcolor{Bittersweet}{ - \int_0^T X_t^{\top}F \bar X_t dt}\right]\\
    &= \mathbb{E}\left[ \frac{1}{2} X_0^{\top}\eta_0 X + X_0^{\top} r_0 +s_0\right]\textcolor{Bittersweet}{-\int_0^T \bar X_t^{\top}F \bar X_t dt}\\
    &= \frac{1}{2} \left( Var(\sqrt{\eta_0} X_0) + \mathbb{E}[\sqrt{\eta_0} X_0]^2 \right) +\bar X_0^{\top} r_0 + s_0\textcolor{Bittersweet}{-\int_0^T \bar X_t^{\top}F \bar X_t dt}.
\end{split}    
\end{align*}
\end{proof}

\begin{proof}[Proof of Theorem~\ref{the:exist_uniq_mfc}]
For the sake of space, we omit the proof of existence and uniqueness, which follows the same steps as in the proof of Theorem~\ref{the:exist_uniq_mfg}.
\end{proof}

\bibliographystyle{siam}

\begin{thebibliography}{10}

\bibitem{achdou_2020}
{\sc Y.~Achdou, C.~Bertucci, J.-M. Lasry, P.-L. Lions, A.~Rostand, and
  J.~Scheinkman}, {\em A class of short-term models for the oil industry
  addressing speculative storage}, arXiv preprint arXiv:2003.11790,  (2020).

\bibitem{huyen_2020}
{\sc R.~A{\"\i}d, M.~Basei, and H.~Pham}, {\em A {M}c{K}ean--{V}lasov approach
  to distributed electricity generation development}, Mathematical Methods of
  Operations Research, 91 (2020), pp.~269--310.

\bibitem{aid_2020}
{\sc R.~A{\"\i}d, R.~Dumitrescu, and P.~Tankov}, {\em The entry and exit game
  in the electricity markets: a mean-field game approach}, arXiv: Optimization
  and Control,  (2020).

\bibitem{alasseur_2020}
{\sc C.~Alasseur, I.~Ben~Taher, and A.~Matoussi}, {\em An extended mean field
  game for storage in smart grids}, Journal of Optimization Theory and
  Applications, 184 (2020), pp.~644--670.

\bibitem{aurell2020optimal}
{\sc A.~Aurell, R.~Carmona, G.~Dayanikli, and M.~Lauriere}, {\em Optimal
  incentives to mitigate epidemics: a {S}tackelberg mean field game approach},
  2020.

\bibitem{malhame_2017}
{\sc O.~Bahn, A.~Haurie, and R.~Malhame}, {\em Limit Game Models for Climate
  Change Negotiations}, 07 2017, pp.~27--47.

\bibitem{bensoussan_yam_2016}
{\sc A.~Bensoussan, M.~H.~M. Chau, and S.~C.~P. Yam}, {\em Mean field games
  with a dominating player}, Applied Mathematics \& Optimization, 74 (2016),
  pp.~91--128.

\bibitem{frehse_2013}
{\sc A.~Bensoussan, J.~Frehse, and P.~Yam}, {\em Mean field games and mean
  field type control theory}, 2013.

\bibitem{caines_huang_malhame_2006}
{\sc P.~E. Caines, M.~Huang, and R.~P. Malhamé}, {\em Large population
  stochastic dynamic games: closed-loop {M}c{K}ean-{V}lasov systems and the
  nash certainty equivalence principle}, Communications in Information and
  Systems, 6 (2006), p.~221–252.

\bibitem{carmona_2018}
{\sc R.~Carmona and F.~Delarue}, {\em FBSDEs and the Solution of MFGs without
  common noise}, Springer International Publishing, Cham, 2018, pp.~215--345.

\bibitem{carmona_ASCONA}
{\sc R.~Carmona and M.~Fehr}, {\em The clean development mechanism and {CER}
  price formation in the carbon emissions markets},  (2010), pp.~341--384.

\bibitem{carmona_SICON}
{\sc R.~Carmona, M.~Fehr, and J.~Hinz}, {\em Optimal stochastic control and
  carbon price formation}, SIAM Journal on Control and Optimization, 48 (2009),
  pp.~2168--2190.

\bibitem{carmona_SIREV}
{\sc R.~Carmona, M.~Fehr, J.~Hinz, and A.~Porchet}, {\em Market design for
  emissions markets trading schemes}, SIAM Review, 52 (2010), pp.~403--452.

\bibitem{wang_2020}
{\sc R.~Carmona and P.~Wang}, {\em Finite-state contract theory with a
  principal and a field of agents}, Management Science,  (2020), p.~First on
  line.

\bibitem{sircar_2017}
{\sc P.~Chan and R.~Sircar}, {\em Fracking, renewables, and mean field games},
  SIAM Review, 59 (2017), pp.~588--615.

\bibitem{djehiche_2018}
{\sc B.~Djehiche, J.~Barreiro-Gomez, and H.~Tembine}, {\em Electricity price
  dynamics in the smart grid : A mean-field-type game perspective}, 2018.

\bibitem{elie_2020}
{\sc R.~{\'E}lie, E.~Hubert, T.~Mastrolia, and D.~Possama{\"\i}}, {\em
  Mean--field moral hazard for optimal energy demand response management},
  Mathematical Finance, 31 (2021), pp.~399--473.

\bibitem{posamai_2019}
{\sc R.~{\'E}lie, T.~Mastrolia, and D.~Possama{\"\i}}, {\em A tale of a
  principal and many, many agents}, Mathematics of Operations Research, 44
  (2019), pp.~440--467.

\bibitem{gueant_2010}
{\sc O.~Gu{\'e}ant, J.-M. Lasry, and P.~L. Lions}, {\em {Mean Field Games and
  oil production}}, in {The Economics of Sustainable Development}, Economica,
  ed., 2010.

\bibitem{hu2019wellposedness}
{\sc K.~Hu}, {\em The wellposedness of path-dependent multidimensional
  forward-backward sde}, 2019.

\bibitem{Jacobson_riccatiboundedness}
{\sc D.~Jacobson}, {\em New conditions for boundedness of the solution of a
  matrix {R}iccati differential equation}, Journal of Differential Equations, 8
  (1970), pp.~258--263.

\bibitem{kucera_2009}
{\sc V.~Ku{\v c}era}, {\em Riccati equations and its solutions}, in Control
  System Advanced Methods, W.~S. Levine, ed., CRC Press, Inc, 2nd~ed., 2009,
  ch.~14.

\bibitem{lasry_lions_2007}
{\sc J.-M. Lasry and P.-L. Lions}, {\em Mean field games}, Japanese Journal of
  Mathematics, 2 (2007), p.~229–260.

\bibitem{AMSnotesLauriere}
{\sc M.~Lauri{\`e}re}, {\em On numerical methods for mean field games and mean
  field type control}, arXiv preprint arXiv:2106.06231,  (2020).

\bibitem{ma_2007}
{\sc J.~Ma and J.~Yong}, {\em Forward-Backward Stochastic Differential
  Equations and their Applications}, Springer Berlin Heidelberg, Berlin,
  Heidelberg, 2007.

\bibitem{salhab2016}
{\sc R.~Salhab, R.~Malham{\'e}, and J.~L. Ny}, {\em A dynamic collective choice
  model with an advertiser}, 2016 IEEE 55th Conference on Decision and Control
  (CDC),  (2016), pp.~6098--6104.

\bibitem{shrivats_2020}
{\sc A.~Shrivats, D.~Firoozi, and S.~Jaimungal}, {\em Optimal generation and
  trading in solar renewable energy certificate (srec) markets}, Applied
  Mathematical Finance, 27 (2020), pp.~131--99.

\bibitem{zhang_2016}
{\sc S.~Zhang, X.~Wang, and A.~Shanain}, {\em Modeling and computation of mean
  field equilibria in producers' game with emission permits trading},
  Communications in Nonlinear Science and Numerical Simulation, 37 (2016),
  pp.~238--248.

\end{thebibliography}

\end{document}